\documentclass{amsproc} \RequirePackage{amssymb}

\copyrightinfo{}{}

\makeatletter

\let\MRhref\@gobble \def\@cite#1#2{{%
 \m@th\upshape\mdseries[{\rmfamily #1}{\if@tempswa, #2\fi}]}}

\theoremstyle{plain} \newtheorem{theorem}{Theorem}[section] 
\newtheorem{proposition}[theorem]{Proposition} 
\newtheorem{lemma}[theorem]{Lemma}

\theoremstyle{definition} 

\theoremstyle{remark} \newtheorem*{remark}{Remark} 
\newtheorem*{examples}{Examples}

\numberwithin{equation}{section}

\newcommand{\union}{\cup} \newcommand{\Union}{\bigcup\limits} 
\newcommand{\inter}{\cap} 

\newcommand{\C}{\mathbb{C}} \newcommand{\N}{\mathbb{N}} 
\newcommand{\R}{\mathbb{R}}  
\newcommand{\Z}{\mathbb{Z}}

\DeclareMathOperator{\codim}{codim}

 \DeclareMathOperator{\Ob}{Ob}

\newcommand{\BDC}{\mathbf{D}^{\mathrm{b}}}

\newcommand{\dsum}[1][]{\mathbin{\oplus_{#1}}}

\DeclareMathOperator{\coker}{coker}

\renewcommand{\to}[1][]{\xrightarrow[#1]{}} 
\newcommand{\from}[1][]{\xleftarrow[#1]{}} 
\newcommand{\isoto}[1][]{\xrightarrow[#1]{\sim}} 
\newcommand{\isofrom}[1][]{\xleftarrow[#1]{\sim}}

\newcommand{\Endo}[1][]{\mathrm{End}_{\raise1.5ex\hbox to.1em{}#1}}

\newcommand{\Hom}[1][]{\mathrm{Hom}_{\raise1.5ex\hbox to.1em{}#1}}

\newcommand{\RHom}[1][]{\mathrm{RHom}_{\raise1.5ex\hbox to.1em{}#1}}

\newcommand{\Ext}[2][]{\mathrm{Ext}_{\raise1.5ex\hbox 
to.1em{}#1}^{#2}}

\newcommand{\THom}[1][]{\mathrm{THom}_{\raise1.5ex\hbox to.1em{}#1}}

\newcommand{\Tens}[1][]{\mathbin{\otimes_{\raise1.5ex\hbox 
to-.1em{}#1}}}

\newcommand{\LTens}[1][]{\mathbin{\otimes_{\raise1.5ex\hbox 
to-.1em{}#1}^{L}}}

\newcommand{\Tor}[2][]{\mathrm{Tor}^{\raise1.5ex\hbox 
to.1em{}#1}_{#2}}

 \def\shb{\mathcal{B}} \def\shc{\mathcal{C}} 
  \def\shf{\mathcal{F}} 
\def\shg{\mathcal{G}} \def\shh{\mathcal{H}}  
 \def\shk{\mathcal{K}}  
\def\shm{\mathcal{M}} \def\shn{\mathcal{N}}  
   
\def\shs{\mathcal{S}}

\newcommand{\sect}{\Gamma} \newcommand{\rsect}{\mathrm{R}\Gamma}

\renewcommand{\hom}[1][]{{\mathcal{H}om}_{\raise1.5ex\hbox 
to.1em{}#1}}

\newcommand{\rhom}[1][]{{R\mathcal{H}om}_{\raise1.5ex\hbox 
to.1em{}#1}}

\newcommand{\ext}[2][]{{\mathcal{E}xt}_{\raise1.5ex\hbox 
to.1em{}#1}^{#2}}

\newcommand{\thom}[1][]{{T\mathcal{H}om}_{\raise1.5ex\hbox 
to.1em{}#1}}

\newcommand{\tens}[1][]{\mathbin{\otimes_{\raise1.5ex\hbox 
to-.1em{}#1}}}

\newcommand{\ltens}[1][]{\mathbin{\otimes_{\raise1.5ex\hbox 
to-.1em{}#1}^{L}}}

\newcommand{\tor}[2][]{{\mathcal{T}or}^{\raise1.5ex\hbox 
to.1em{}#1}_{#2}}

\newcommand{\wtens}{\mathbin{\mathop{\otimes}\limits^{{}_{\mathrm{w}}}}}

\DeclareMathOperator{\supp}{supp}

\newcommand{\oim}[1]{{#1}_*} 

 \newcommand{\reim}[1]{{R#1}_!}

\newcommand{\opb}[1]{#1^{-1}}

\DeclareMathOperator{\ori}{or}

\newcommand{\GHom}[1][]{\mathrm{GHom}_{\raise1.5ex\hbox to.1em{}#1}}

\newcommand{\GExt}[2][]{\mathrm{GExt}_{\raise1.5ex\hbox 
to.1em{}#1}^{#2}}

\newcommand{\FHom}[1][]{\mathrm{FHom}_{\raise1.5ex\hbox to.1em{}#1}}

\newcommand{\ghom}[1][]{{\mathcal{GH}om}_{\raise1.5ex\hbox 
to.1em{}#1}}

\newcommand{\gext}[2][]{{\mathcal{GE}xt}_{\raise1.5ex\hbox 
to.1em{}#1}^{#2}}

\newcommand{\fhom}[1][]{{\mathcal{FH}om}_{\raise1.5ex\hbox 
to.1em{}#1}}

\newcommand{\tenstop}[1][]{\mathbin{\hat{\otimes}_{\raise1.5ex\hbox 
to-.1em{}#1}}}

\newcommand{\homtop}[1][]{\mathcal{L}_{\raise1.5ex\hbox to.1em{}#1}}

\newcommand{\Homtop}[1][]{\mathrm{L}_{\raise1.5ex\hbox to.1em{}#1}}

\newcommand{\D}{\mathcal{D}}

\renewcommand{\O}{\mathcal{O}}

\newcommand{\B}{\mathcal{B}}

\newcommand{\Db}{\mathcal{D}b}

\DeclareMathOperator{\chv}{char}

\newcommand{\dsol}{\mathop{\mathcal{S}ol}\nolimits} 
\newcommand{\DSol}{\mathop{\mathrm{Sol}}\nolimits}

\newcommand{\detens}{\mathbin{\underline{\boxtimes}}}

\def\absdoim#1{\underline{#1}_*} 
\def\reldoim[#1]#2{\underline{#2}_{|{#1}*}} \def\doim{\@ifnextchar 
[{\reldoim}{\absdoim}}

\def\absdeim#1{\underline{#1}_*} 
\def\reldeim[#1]#2{\underline{#2}_{|{#1}*}} \def\deim{\@ifnextchar 
[{\reldeim}{\absdeim}}

\def\absdopb#1{\underline{#1}^{-1}} 
\def\reldopb[#1]#2{\underline{#2}_{|{#1}}^{-1}} \def\dopb{\@ifnextchar 
[{\reldopb}{\absdopb}}

\def\absboim#1{\underline{\underline{#1}}_*} 
\def\relboim[#1]#2{\underline{\underline{#2}}_{|{#1}*}} 
\def\boim{\@ifnextchar [{\relboim}{\absboim}}

\def\absbeim#1{\underline{\underline{#1}}_*} 
\def\relbeim[#1]#2{\underline{\underline{#2}}_{|{#1}*}} 
\def\beim{\@ifnextchar [{\relbeim}{\absbeim}}

\def\absbopb#1{\underline{\underline{#1}}^*} 
\def\relbopb[#1]#2{\underline{\underline{#2}}_{|{#1}}^*} 
\def\bopb{\@ifnextchar [{\relbopb}{\absbopb}}

\newcommand{\ddual}{\underline{D}}

\newcommand{\coh}{\mathrm{coh}}

\newcommand{\rcons}{{\R\mathrm{-c}}}

\makeatletter

\newcommand{\shift}[1]{{\scriptstyle{[#1]}}} 
\newcommand{\twist}[1]{{\scriptstyle{(#1)}}} 
\newcommand{\abstwist}[1]{{\scriptstyle{\{#1\}}}}

\renewcommand{\hom}[1][]{\operatorname{hom}\nolimits_{#1}} 
\renewcommand{\rhom}{\mathrm{R}\!\hom} \renewcommand{\thom}{t\!\hom}

\renewcommand{\deim}[1]{{D #1}_{!}} \renewcommand{\dopb}[1]{D #1^{*}} 
\newcommand{\dtens}{\mathbin{\mathop{\otimes}\limits^{{}_{D}}}} 
\renewcommand{\detens}{\mathbin{\mathop{\boxtimes}\limits^{{}_{D}}}} 
\newcommand{\dcirc}{\mathbin{\mathop{\circ}\limits^{{}_{D}}}} 
\renewcommand{\ddual}[1]{{#1}^{\vee}}

\newcommand{\OX}{\O_X}  
\newcommand{\OXY}{\O_{X\times \Xi}} \newcommand{\DX}{\D_X} 
\newcommand{\DY}{\D_\Xi}

\newcommand{\GL}{\mathop{\mathrm{GL}}}

\newcommand{\CV}{\mathbb{V}} \newcommand{\RV}{\mathsf{V}} 
\newcommand{\MV}{W}

\renewcommand{\AA}{\mathbb{A}} \renewcommand{\aa}{\AA^*} 
\newcommand{\PP}{\mathbb{P}} \newcommand{\bb}{\PP^*} 
\newcommand{\GG}{\mathbb{G}} \newcommand{\MM}{\mathbb{M}} 
\newcommand{\FF}{\mathbb{F}} 

\newcommand{\OP}{\O_{\PP}} \renewcommand{\Ob}{\O_{\bb}} 
\newcommand{\OG}{\O_{\GG}} \newcommand{\DP}{\D_{\PP}} 
\renewcommand{\Db}{\D_{\bb}} \newcommand{\DG}{\D_{\GG}} 
\newcommand{\OPG}{\O_{\PP\times\GG}}

\newcommand{\UP}{\mathsf{Q}} \newcommand{\UG}{\mathsf{M}}

\newcommand{\RP}{\mathsf{P}} \newcommand{\RG}{\mathsf{G}} 
\newcommand{\Rb}{\RP^*} \newcommand{\RA}{\mathsf{A}}

\newcommand{\RGA}{\RG_{\RA}} \newcommand{\RH}{\mathsf{H}} 
\newcommand{\HH}{\mathbb{H}} \newcommand{\Rh}{\RH^*} 
\newcommand{\hh}{\HH^*} \newcommand{\Rbo}{\Rb_{\circ}} 
\newcommand{\bbo}{\PP^*_{\circ}} \newcommand{\qRb}{\tilde\Rb} 
\newcommand{\qRbo}{\tilde\Rbo} \newcommand{\xio}{h}

\newcommand{\BF}{\shb_{\FF}}

\renewcommand{\dsol}{\operatorname{sol}\nolimits}

\title{Sheaves and $\D$-modules in integral geometry%
\footnotetext{{\bf Appeared in:} 
Analysis,
  geometry, number theory: the mathematics of Leon Ehrenpreis (Philadelphia,
  PA, 1998), Amer. Math. Soc., Providence, RI, 2000, Contemp. Math., 251,
  pp.~141--161.}
}

\author{Andrea D'Agnolo}

\begin{document}

\begin{abstract}
    Integral geometry deals with those integral transforms which
    associate to ``functions'' on a manifold $X$ their integrals along
    submanifolds parameterized by another manifold $\Xi$.  Basic
    problems in this context are range characterization---where
    systems of PDE appear---and inversion formulae.  As we pointed out
    in a series of joint papers with Pierre Schapira, the language of
    sheaves and $\D$-modules provides both a natural framework and
    powerful tools for the study of such problems.  In particular, it
    provides a general adjunction formula which is a sort of
    archetypical theorem in integral geometry.  Focusing on range
    characterization, we illustrate this approach with a discussion of
    the Radon transform, in some of its manifold manifestations.
\end{abstract}
    
\maketitle

\section*{Introduction}

According to I.~M.~Gelfand, integral geometry is the study of those
integral transforms $\varphi(x)\mapsto\psi(\xi)=\int\varphi(x)
k(x,\xi)$ which are associated to a geometric correspondence $S\subset
X\times \Xi$ between two manifolds $X$ and $\Xi$.  In other words,
given the family of subvarieties $\widehat \xi=\{x\colon (x,\xi)\in
S\}\subset X$, parameterized by the manifold $\Xi\owns\xi$, one is
interested in those transforms whose integral kernel $k$ is a characteristic
class of the incidence relation $S$.  Typical examples are the various
instances of the Radon transform.  For the real affine Radon
transform, $\widehat\xi$ is a $p$-plane in an affine space $X$, and
$\psi(\xi)=\int_{\widehat\xi}\varphi(x)$ is the integral of a rapidly
decreasing $C^\infty$-function with respect to the Euclidean measure. 
For the complex projective Radon transform, $\widehat\xi$ is a $p$-plane
in a complex projective space $X$, $\varphi$ is a cohomology class of
a holomorphic line bundle, and the integral $\int_{\widehat\xi}$ has to be
understood as a Leray-Grothendieck residue.  The conformal Radon
transform may be viewed as a boundary value of the complex projective
transform and, in this case, the natural function space to consider is
that of hyperfunctions.  We begin by recalling in the first section
some classical results on range characterization for these Radon
transforms, where systems of partial differential equations appear.

In the second and main part of this survey, we present a general
framework for integral geometry, based on the theory of sheaves and
$\D$-modules, which we proposed in a series of joint papers with
Pierre Schapira.  This approach makes apparent the common pattern
underlying the various instances of the Radon transform.  Let us
discuss it in some detail.

It is a central idea of M.~Sato to obtain hyperfunctions on a real
analytic manifold $M$ as boundary value of holomorphic functions in a
complex neighborhood $X\supset M$.  More generally, one may obtain
generalized functions by coupling holomorphic functions with a
constructible sheaf $F$ on $X$.  Following Kashiwara-Schapira, we
consider four examples of such generalized functions, that we denote
by $\shc^{\omega}(F)$, $\shc^{\infty}(F)$, $\shc^{-\infty}(F)$,
$\shc^{-\omega}(F)$.  If $F$ is the constant sheaf along $M$, these
are the sheaves of analytic functions, $C^\infty$-functions,
distributions and hyperfunctions, respectively.  Let thus $X$ and
$\Xi$ be complex manifolds.  An integral transform $\int\varphi(x)
k(x,\xi)$ consists of three operations. Namely, 
a pull-back from $X$ to
$X\times \Xi$, a product with an integral kernel, and a push-forward to $\Xi$. 
The Grothendieck formalism of operations makes sense in the categories
of sheaves and $\D$-modules.  Let us denote here for short by
$\widehat\shm$ the transform of a coherent $\D$-module $\shm$ on $X$,
and by $\widehat G$ the transform of a constructible sheaf $G$ on
$\Xi$.  A general result in this framework asserts that, under mild
hypotheses, solutions of $\shm$ with values in
$\shc^{\natural}(\widehat G)$ are isomorphic to solutions of
$\widehat\shm$ with values in $\shc^{\natural}(G)$, for
$\natural=\pm\infty,\ \pm\omega$.  This shows that, in order to deal
with concrete examples, one has to address three problems of an
independent nature
\begin{itemize}
    \item [(T)] the computation of $\widehat G$,

    \item [(A)] the computation of $\widehat\shm$,

    \item [(Q)] the explicit description of the isomorphism between
    solutions.
\end{itemize}
Problem (T) is the easiest one, since it is of a topological nature. 
Problem (A) is of an analytical nature, and we will explain how the
theory of microlocal operators helps in dealing with it.  In
particular, under certain assumptions on the microlocal geometry
associated to the correspondence (which are satisfied by the Radon
transform), the functor $\shm\mapsto\widehat\shm$ establishes an
equivalence of categories between coherent $\D$-modules on $X$ and
systems with regular singularities along an involutive submanifold of
the cotangent bundle to $\Xi$, determined by the correspondence. 
Finally, (Q) is a kind of quantization problem, which consists in
exhibiting the integral kernel of the transform.  This is obtained as
boundary value of the meromorphic kernel associated to the $\D$-module
transform, which may be viewed as a complex analogue of Lagrangian
distributions.

In the third and last section of our survey, we exemplify the above
approach by giving proofs of the results on range characterization for
Radon transforms.  These are taken from joint works with P.~Schapira,
with C.~Marastoni, or by the present author alone.  Although we
restricted our attention to classical results, it should already be
evident that the tools of sheaf and $\D$-module theory are quite
powerful in this context.  In particular, they allow one to treat
higher cohomology groups, obtaining results which are hardly accessible
by classical means.  As an example, we get a geometrical
interpretation of the Cavalieri condition, which characterizes the
image of the real affine Radon transform.

Let us make some comments on related results.  Independently, several
people proposed approaches to integral geometry which are more or less
close to ours.  A.~B.~Goncharov gave a $\D$-module interpretation of
Gelfand's $\kappa$-form, used to describe inversion formulae. 
Kashiwara-Schmid announced an adjunction formula similar to ours,
taking into account topology and group actions.  In an equivariant
setting, the computation of the $\D$-module transform for
correspondences between generalized flag manifolds appear in recent
works by T.~Oshima and T.~Tanisaki.  The computations in
Baston-Eastwood should also be useful in this context.

Concerning the exposition that follows, to keep it lighter we decided
to postpone the bibliographical comments to the end of each section. 
As a consequence, a lack of reference does not imply that a result
should be attributed to the present author.

\section{Radon transform(s)}
\label{se:Statements}

Here, we collect some more or less classical results on range
characterization for various instances of the Radon $p$-plane
transform.  This should be considered as a motivation for the theory
we present in the next section.

\subsection{Real affine case}
\label{ss:realaffine}

The spaces
\begin{equation}
    \begin{cases}
	\RA &\text{a real $n$-dimensional affine space,} \\
	\RGA &\text{the family of affine $p$-planes in $\RA$,}
    \end{cases}
    \label{eq:geoRA}
\end{equation}
are related by the correspondence
$\RGA\owns\xi\multimap\widehat\xi\subset\RA$, associating to a plane
the set of its points.  Denote by $\shs(\RA)$ the Schwartz space of
rapidly decreasing $C^\infty$-functions on $\RA$.  The classical
affine Radon transform consists in associating to a function $\varphi$
in $\shs(\RA)$ its integrals along the family of $p$-planes
\begin{equation}
    R_{\RA}\colon\shs(\RA)\owns\varphi\mapsto 
    \psi(\xi)=\int_{\widehat\xi}\varphi.
    \label{eq:Raffine}
\end{equation}
In order to explain what measure is used in the above integral, let us 
choose a system of coordinates $(t)=(t_{1},\dots,t_{n})$ in 
$\RA\simeq\R^n$.  Any affine $p$-plane $\xi\in\RGA$ can be described 
by a system of equations
$$
\xi\colon \langle t,\tau_{i}\rangle + \sigma_{i} =0, \qquad 
i=1,\dots,n-p,
$$
where $\tau_{i}=(\tau_{i}^1,\dots,\tau_{i}^n)$ are linearly 
independent vectors in $(\R^n)^*$, and $\sigma_{i}\in\R$.  Consider
$$
\psi(\sigma,\tau) = \int \varphi(t) \delta(\langle t,\tau_{1} \rangle 
+ \sigma_{1}) \cdots \delta(\langle t,\tau_{n-p} \rangle + 
\sigma_{n-p}) \, dt_{1}\cdots dt_{n},
$$
where $\sigma\in\R^{n-p}$ is the row vector 
$\sigma=(\sigma_{1},\dots,\sigma_{n-p})$, 
$\tau=(\tau_{1},\dots,\tau_{n-p})$ is an $n\times (n-p)$ matrix, and 
$\delta$ denotes the Dirac delta function.  Since $\psi$ satisfies the 
homogeneity condition
$$
\psi(\sigma\mu,\tau\mu)=|\det\mu|^{-1} \psi(\sigma,\tau) \quad 
\forall\mu\in\GL(n-p,\R),
$$
it defines a section of a line bundle over $\RGA$, that we denote by 
$\shc^\infty_{\RGA}\abstwist{-1}$.  The projection
$$
(\xi\colon \langle t,\tau_{i}\rangle + \sigma_{i} =0) \mapsto 
(\xi'\colon \langle t,\tau_{i}\rangle = 0)
$$
makes $\RGA$ into an $(n-p)$-dimensional vector bundle 
$q\colon\RGA\to\RG'$ over the compact Grassmannian of vector 
$p$-planes in $\R^n$.  One says that a global section of 
$\shc^\infty_{\RGA}\abstwist{-1}$ is rapidly decreasing, if it is 
rapidly decreasing along the fibers of $q$.

\begin{theorem}
    \label{th:Rreal}
    The real affine Radon transform
    \begin{align*}
	R_{\RA}\colon\shs(\RA) &\to 
	\sect(\RGA;\shc^\infty_{\RGA}\abstwist{-1}) \\
	\varphi(t) &\mapsto \psi(\sigma,\tau)
    \end{align*}
    is injective.  Concerning its range
    \begin{itemize}
	\item [(i)] If $p<n-1$, then $\psi$ belongs to the image of 
	$R_{\RA}$ if and only if it is rapidly decreasing and 
	satisfies the John ultrahyperbolic system
	\begin{equation}
	    \label{eq:John}
	\left(\frac {\partial^2} {\partial\xi_{i}^{j} 
	\partial\xi_{i'}^{j'}} - \frac {\partial^2} 
	{\partial\xi_{i}^{j'} \partial\xi_{i'}^{j}}\right) \psi(\xi) = 
	0, \qquad
	\begin{array}{l}
	    j,j'=0,\dots,n, \\
	    i,i'=1,\dots,n-p,
	\end{array}
	\end{equation}
	where $(\xi)=(\xi_{1},\dots,\xi_{n-p})$, $\xi_{i} = 
	(\xi_{i}^0,\dots\xi_{i}^n) = 
	(\sigma_{i},\tau_{i}^1,\dots,\tau_{i}^n)$.
    
	\item [(ii)] If $p=n-1$, then $\psi$ belongs to the image of 
	$R_{\RA}$ if and only if it is rapidly decreasing and 
	satisfies the Cavalieri condition $$
	\int_{-\infty}^{+\infty}\psi(\sigma,\tau)\sigma^m\,d\sigma 
	\quad\text{is, for any $m\in\N$, a polynomial in $\tau$ of 
	degree $\leq m$.} $$
	(Since $n-p=1$, we have here $\sigma=\sigma_{1}$ and 
	$\tau=\tau_{1}$.)
    \end{itemize}
\end{theorem}

Heuristically, the system of PDE in~(i) compensates for the difference 
of dimensions between $\dim\RA=n$ and $\dim\RGA=(p+1)(n-p)$.  Also for 
$p=n-1$ some conditions should have been expected.  In fact, 
$\shs(\RA)$ is the space of functions which are rapidly decreasing in 
the $n-1$ directions of $\RA$ going to infinity, while rapid decrease 
in $\RGA$ concerns only the $1$-dimensional fiber of 
$q\colon\RGA\to\RG'$.  The Cavalieri condition---also known as moment 
condition---was classically obtained using Fourier inversion formula.  
Instead, using our formalism, we will obtain this condition in a 
purely geometric way.

\subsection{Complex projective case}
\label{ss:complexproj}

Let us consider
\begin{equation}
    \begin{cases}
	\CV &\text{a complex vector space of dimension $n+1$,} \\
	\PP &\text{the projective space of complex vector lines in 
	$\CV$,} \\
	\GG &\text{the Grassmannian of projective $p$-planes in 
	$\PP$.}
    \end{cases}
    \label{eq:geoC}
\end{equation}
As for the real case, there is a correspondence 
$\GG\owns\zeta\multimap\widehat\zeta\subset\PP$, associating to a plane 
the set of its points.  Reciprocally, denote by $\widehat z\subset\GG$ the 
set of $p$-planes containing $z\in\PP$.  We also set $\widehat 
U=\Union\nolimits_{\zeta\in U}\widehat\zeta$ for $U\subset\GG$.

If $[z]=[z_{0},\dots,z_{n}]$ are homogeneous coordinates in $\PP$, a 
projective $p$-plane is described by the system of equations
\begin{equation}
    \zeta\colon \langle z,\zeta_{i}\rangle = 0, \qquad i=1,\dots,n-p,
    \label{eq:zeta}
\end{equation}
where $\zeta_{i}$ are linearly independent vectors in $\CV^*$.

For $m\in\Z$, let us denote by $\OP\twist{m}$ the holomorphic line 
bundle whose sections $\varphi$, written in homogeneous coordinates, 
satisfy the homogeneity condition
$$
\varphi(\lambda z)=\lambda^m \varphi(z) \quad 
\forall\lambda\in\C^\times=\GL(1,\C).
$$
The Leray form
$$
\omega(z) = \sum_{j=0}^{n} (-1)^j z_j dz_0\wedge\cdots\wedge dz_{j-1} 
\wedge dz_{j+1} \wedge\cdots\wedge dz_n
$$
is $(n+1)$-homogeneous, i.e.\ is a global section of 
$\Omega_{\PP}\twist{n+1}=\Omega_{\PP}\tens[\OP]\OP\twist{n+1}$, where 
$\Omega_{\PP}$ denotes the sheaf of holomorphic forms of maximal 
degree.

Let $U\subset\GG$ be an open subset such that $\widehat z\inter U$ is 
connected for any $z\in\widehat U$.  Consider the integral transform
\begin{equation}
    \label{eq:kRC}
    R_{\PP}\colon H^p(\widehat U;\OP\twist{-p-1})\owns\varphi\mapsto 
    \psi(\zeta) = \left(\!\frac 1 {2\pi i}\!\right)^{\!\!\!  n-p} 
    \!\!\int\varphi(z)\frac{\omega(z)}{\langle z,\zeta_{1}\rangle 
    \cdots \langle z,\zeta_{n-p}\rangle},
\end{equation}
where the integral sign stands for a Leray-Grothendieck residue.  Note 
that the integrand is a well defined function in the $\PP$ variable, 
since it is homogeneous of degree zero in $z$.  Moreover, $\psi(\zeta) 
\in \sect(U;\OG\twist{-1})$, where $\OG\twist{-1}$ is the line bundle 
whose sections satisfy the homogeneity condition
$$
\psi(\zeta\mu)=(\det\mu)^{-1} \psi(\zeta) \quad 
\forall\mu\in\GL(n-p,\C).
$$
Recall that, by the Cauchy formula, the Dirac delta function 
$\delta(u)$ in $\R$ is the boundary value of $(2\pi i \, w)^{-1}$, for 
$w=u+iv\in\C$.  Thus, the integral kernel of $R_{\PP}$ is the complex 
analog of the delta function along the incidence relation 
$z\in\widehat\zeta$.  In this sense $R_{\PP}$ is indeed a complex Radon 
transform $\varphi\mapsto\int_{\widehat\zeta}\varphi$.

Following Penrose, if $p=1$ and $n=3$ then $\GG=\MM$ is a conformal 
compactification of the linear complexified Minkowski space.  
Maxwell's wave equation thus extends to an equivariant differential 
operator on $\MM$.  For arbitrary $p<n-1$, we denote by $\square$ its 
higher dimensional analog in $\GG$
\begin{equation}
    \label{eq:wave}
    \square\colon\OG\twist{-1}\to\shh,
\end{equation}
where $\shh$ is a homogeneous vector bundle that we do not need to 
describe here.  Let us just point out that $\square$ is a complex 
version of the John ultrahyperbolic system.  In fact, 
$\square\psi(\zeta)$ is expressed by the left hand side term of 
\eqref{eq:John}, replacing $(\xi_{i}^j)$ with the dual Stiefel 
coordinates of $\GG$, 
$(\zeta_{i}^j)\in (\CV^*)^{n-p}$.

One says that $U\subset\GG$ is elementary if for any $z\in\widehat U$
the slice $\widehat z\inter U$ has trivial reduced cohomology up to
degree $p$.  This means that $\widehat z\inter U$ should be connected,
and its Betti numbers $b^j$ should vanish for $1\leq j\leq p$.

\begin{theorem}
    \label{th:Rcomplex}
    Let $U$ be an elementary open subset of $\GG$, and consider the 
    complex Radon transform $$
    R_{\PP}\colon H^p(\widehat U;\OP\twist{-p-1}) \to 
    \sect(U;\OG\twist{-1}).  $$
    \begin{itemize}
	\item [(i)] If $p<n-1$, then $R_{\PP}$ is an isomorphism onto 
	the space of sections $\psi$ satisfying $\square\psi=0$.
    
	\item [(ii)] If $p=n-1$, then $\GG=\bb$ is a dual projective 
	space, and $R_{\PP}$ is an isomorphism.
    \end{itemize}
\end{theorem}

\subsection{Real conformal case}
\label{ss:realconformal}

For simplicity sake, let us restrict our attention to a generalization 
of the Penrose case, where $n+1 = 2(p+1)$.  Consider the following 
geometrical situation
\begin{equation}
    \begin{cases}
	\Phi &\text{a Hermitian form on $\CV\simeq\C^{2(p+1)}$ of 
	signature $(p+1,p+1)$,} \\
	\UP\subset\PP &\text{the set of null vectors $z$, with 
	$\Phi\vert_{z}=0$,} \\
	\UG\subset\GG &\text{the set of null planes $\zeta$, with 
	$\Phi\vert_{\zeta}=0$.}
    \end{cases}
    \label{eq:geoRU}
\end{equation}
One checks that $\UP$ is a real hypersurface of $\PP$ whose Levi form 
has rank $2p$ and signature $(p,p)$, and $\UG$ is a totally real 
analytic submanifold of $\GG$.  Let us denote by $\shc^\omega_{\UG}$ 
and $\shc^{-\omega}_{\UG}$ the sheaves of real analytic functions and 
Sato hyperfunctions on $\UG$, respectively.  Since the higher 
dimensional Maxwell wave operator \eqref{eq:wave} is hyperbolic for 
$\UG$, it is natural to consider its hyperfunction solutions.  For 
$\natural=\pm\omega$, denote by $\ker(\square;\shc^\natural_{\UG})$ 
the subsheaf of $\shc^\natural_{\UG}\tens[\OG]\OG\twist{-1}$ whose 
sections $\psi$ satisfy $\square\psi=0$.

One has the following boundary value version of 
Theorem~\ref{th:Rcomplex}.

\begin{theorem}
    \label{th:hyp}
    There is a commutative diagram $$
    \begin{array}{ccc}
	H^p(\UP;\OP\twist{-p-1}) &\isoto& 
	\sect(\UG;\ker(\square;\shc^{\omega}_{\UG})) \\
	\downarrow && \downarrow \\
	H^{p+1}_{\UP}(\PP;\OP\twist{-p-1}) &\isoto& 
	\sect(\UG;\ker(\square;\shc^{-\omega}_{\UG}))
    \end{array}
    $$
    where the horizontal arrows are induced by $R_{\PP}$, the vertical 
    arrows are the natural maps, and $H^j_{\UP}(\cdot)$ denotes the 
    $j$-th cohomology group with support on $\UP$.
\end{theorem}

\subsection*{Notes}

We give here some references to papers where the Radon transform is dealt
with using a classical approach, i.e.~without using the tools of sheaf
and $\D$-module theory.

\medskip\noindent
{\bf \S\ref{ss:realaffine}} The study of the hyperplane Radon 
transform in affine $3$-space goes back to Radon 
himself~\cite{Radon17}.  Line integrals were considered by 
John~\cite{John38}, where the ultrahyperbolic equation first appeared.  
Theorem~\ref{th:Rreal} may be found, for example, 
in~\cite{Gelfand-Gindikin-Graev82} or \cite{Helgason84}.  See 
also~\cite{Grinberg85}, \cite{Gonzalez91}, and \cite{Kakehi97} for 
related results.

\medskip\noindent {\bf \S\ref{ss:complexproj}} Part~(ii) of
Theorem~\ref{th:Rcomplex} is due to Martineau~\cite{Martineau67},
while part~(i) for $p=1$, $n=3$ is known as the Penrose transform, and
was obtained in~\cite{Eastwood-Penrose-Wells81}.  Correspondences
between generalized flag manifolds are discussed
in~\cite{Baston-Eastwood89}.  Theorem~\ref{th:Rcomplex} is also
considered by~\cite{Gindikin-Henkin78}, \cite{Henkin-Polyakov87}, in
the case where $\widehat U$ is $p$-linearly concave.  The case where
$n+1 = 2(p+1)$ and $U=\{\zeta\colon \Phi\vert_\zeta\gg 0\}$ is
discussed in~\cite{Sekiguchi96}, and generalized to other
correspondences between Grassmannians, using representation
theoretical arguments.

\medskip\noindent
{\bf \S\ref{ss:realconformal}} Theorem~\ref{th:hyp} for $p=1$ is due 
to~\cite{Wells81}.  Note also that the problem of extending to $\UG$ 
sections of $\ker(\square;\shc^{-\omega}_{\UG})$ defined on the affine 
Minkowski space is discussed in~\cite{Bailey-Ehrenpreis-Wells82}.

\section{Sheaves and $\D$-modules for integral geometry}

There is a common pattern underlying the examples we discussed in the 
previous section, which the language of sheaves and $\D$-modules will 
make apparent.  As a motivation for our choice of framework, let us 
state an informal paradigm for integral transforms.

Let $X$ be a manifold, denote by $\shf$ a space of ``functions'' on $X$, 
let $\shm$ represent some system of PDE (possibly void, since the 
absence of differential equations corresponds to the system $0u=0$), 
and denote by $\DSol(\shm,\shf)$ the space of $\shf$-valued solutions 
to $\shm$.  Similarly, consider the space $\DSol(\shn,\shg)$ on 
another manifold $\Xi$.  An integral transform from $X$ to $\Xi$ is a 
map
\begin{eqnarray}
	\DSol(\shm,\shf) &\to& \DSol(\shn,\shg) \label{eq:mfng}\\
	\varphi(x) &\mapsto& \psi(\xi)= \int\varphi(x)\cdot k(x,\xi), 
	\nonumber
\end{eqnarray}
consisting of the composition of three operations
\begin{itemize}
	\item [(i)] pull back $\varphi$ to $X\times \Xi$,

	\item [(ii)] take the product with an integral kernel $k$,

	\item [(iii)] push it forward to $\Xi$.
\end{itemize}

To make some sense out of this, we have to decide what we mean by 
space of functions, how we represent systems of PDE, and how we 
perform the three steps above in this setting.

\subsection{Generalized functions}
\label{ss:GeneralizedFcts}

Let $X$ be a complex manifold, and denote by $\OX$ its structural 
sheaf of holomorphic functions.  Starting from $\OX$, some spaces of 
functions are naturally associated to the datum of a locally closed 
subanalytic subset.

\begin{examples}
$\bullet$ To an open subset $U\subset X$ is associated the space 
$\sect(U;\OX)$ of holomorphic functions on $U$. More generally, one 
may consider the 
cohomology groups $H^j(U;\OX)$.

$\bullet$ To a real analytic submanifold $M$, of which $X$ is a 
complexification, are 
associated the sheaves
\begin{equation}
    \shc^{\omega}_{M} \subset \shc^{\infty}_{M} \subset 
    \shc^{-\infty}_{M} \subset \shc^{-\omega}_{M}
    \label{eq:Mspaces}
\end{equation}
of real analytic functions, $C^\infty$-functions, Schwartz 
distributions and Sato hyperfunctions, respectively.  Recall that, 
locally, hyperfunctions are represented by finite sums $\sum_{j} 
b(F_{j})$, where $b(F_{j})$ denotes the boundary value of a 
holomorphic function $F_{j}$ defined on an open wedge of $X$ with edge 
$M$.  For example, hyperfunctions on the real line are equivalence 
classes in $\sect(\C\setminus\R;\OX)/\sect(\C;\OX)$.  More precisely, 
$\shc^{-\omega}_{M}=H^{\dim X}_{M}(\OX)\tens\ori_{M/X}$ is obtained 
from $\OX$ by considering cohomology with support in $M$, twisted by 
the relative orientation sheaf $\ori_{M/X}$ (since $\ori_{M/X}$ is 
locally trivial, the twisting operation is locally void).
\end{examples}

\noindent The set of spaces as above, obtained from $\OX$ via the 
datum of a subanalytic subset, is not stable by direct images.  For 
example

\medskip $\bullet$ Let $X=\C\setminus\{0\}$ with holomorphic
coordinate $z$, and define $p\colon X\to X$ by $p(z)=z^2$.  Then
$\sqrt z\in\oim p\OX$ is described by a local system rather than by a
subset.  In fact, one has $\oim p\OX\simeq\hom(\oim p\C_{X},\OX)$, and
$\sqrt z\in\hom(L,\OX)$, where $\C_{X}$ denotes the constant sheaf on
$X$, and $L$ is defined by $\oim p\C_{X}\simeq\C_{X}\dsum L$. 
\medskip

Instead of subanalytic subsets, we thus consider the bounded derived 
category $\BDC_{\rcons}(\C_{X})$ of $\R$-constructible sheaves.  
Roughly speaking, this is the smallest full subcategory of the bounded 
derived category of sheaves such that: (i) it contains the skyscraper 
sheaves $\C_{U}$ along open subanalytic subsets $U\subset X$, (ii) the 
assignment $X\mapsto \BDC_{\rcons}(\C_{X})$ is stable by exterior 
tensor products, inverse images, and proper direct images.

With these notations, some of the above examples are expressed as
\begin{align*}
    H^j(U;\OX) &= H^j\RHom(\C_{U},\OX), & \sqrt z &\in \hom(L,\OX), \\
    \shc^{\omega}_{M} &= 
    \C_{M}\tens\OX, & \shc^{-\omega}_{M} &= 
    \rhom(\C_{M}',\OX),
\end{align*}
where we denote by $F'=\rhom(F,\C_{X})$ the dual of
$F\in\BDC_{\rcons}(\C_{X})$.

It is also possible to obtain in a similar way the sheaves of 
$C^\infty$-functions and Schwartz distributions, replacing $\tens$ and 
$\rhom$ by the functors $\smash\wtens$ and $\thom$ of formal and 
tempered cohomology.  Let us briefly recall their construction.  
First, using a result of Lojasiewicz, one proves that there exist 
exact functors $w$ and $t$ in $\BDC_{\rcons}(\C_{X})$, characterized 
by the following requirements.  If $S$ is a closed subanalytic subset 
of $X$, then $w(\C_{X\setminus S})$ is the ideal of $\shc^\infty_X$ of 
functions vanishing to infinite order on $S$, and $t(\C_{S})$ is the 
subsheaf of $\shc^{-\infty}_{X}$ whose sections have support contained 
in $S$.  Then, for $F\in\BDC_{\rcons}(\C_{X})$ one defines 
$F\smash\wtens\OX$ and $\thom(F,\OX)$ as the Dolbeault complexes with 
coefficients in $w(F)$ and $t(F)$, respectively.

To $F\in\BDC_{\rcons}(\C_{X})$ we thus associate the spaces of 
``generalized functions''
$$
\shc^{\natural}(F) =
\begin{cases}
    F\tens\OX & \natural=\omega, \\
    F\wtens\OX & \natural=\infty, \\
    \thom(F',\OX) & \natural=-\infty, \\
    \rhom(F',\OX) & \natural=-\omega.
\end{cases}
$$
In particular, \eqref{eq:Mspaces} is obtained from the natural sequence
$$
    \shc^{\omega}(F) \to \shc^{\infty}(F) \to \shc^{-\infty}(F) \to 
    \shc^{-\omega}(F)
$$
by taking $F=\C_{M}$.

\subsection{Systems of PDE}
\label{ss:PDE}

Denote by $\DX$ the sheaf of linear differential operators with 
holomorphic coefficients in $X$.  By definition, a (left) $\DX$-module 
$\shm$ is coherent if it is locally represented by a system of PDE, 
i.e.~if locally there exists an exact sequence
$$
(\DX)^{N_{1}}\to[\cdot P](\DX)^{N_{0}}\to\shm\to 0,
$$
where $P=(P_{ij})$ is an $N_{1}\times N_{0}$ matrix with elements in 
$\DX$, acting by multiplication to the right $Q\mapsto QP$.  We denote 
by $\BDC_{\coh}(\DX)$ the bounded derived category of $\DX$-modules 
with coherent cohomology groups.  To $\shm\in\BDC_{\coh}(\DX)$ is 
associated its characteristic variety $\chv(\shm)$, a closed complex 
subvariety of the cotangent bundle $T^*X$, which is involutive for the 
natural symplectic structure.

If $\shf$ is another $\DX$-module (not necessarily coherent), the 
solutions of $\shm$ with values in $\shf$ are obtained by
$$
\dsol(\shm,\shf)=\rhom[\DX](\shm,\shf),
$$
and we will also consider the space of global sections
$$
\DSol(\shm,\shf)=\rsect(X;\dsol(\shm,\shf)).
$$

\begin{examples}
$\bullet$ To a holomorphic vector bundle $\shh$ one associates the 
$\DX$-module $\D\shh^*=\DX\tens[\OX]\shh^*$, where 
$\shh^*=\hom[\OX](\shh,\OX)$.  Since $\D\shh^*$ is locally represented 
by the trivial system $0u=0$, its characteristic variety is the whole 
space $T^*X$.  We can recover the $\C$-module underlying $\shh$ by 
taking holomorphic solutions to $\D\shh^*$, namely 
$\shh\simeq\dsol(\D\shh^*,\OX)$.

$\bullet$ Next, consider the cyclic $\DX$-module $\DX^P=\DX/\DX P$ 
associated to a single differential operator $P\in\DX$.  Its 
characteristic variety is the hypersurface defined by the zero locus 
of the principal symbol of $P$.  Concerning its solutions, 
$H^0\dsol(\DX^P,\shf)=\{\varphi\in\shf\colon P\varphi=0\}$ describes 
the kernel of the operator $P$ acting on $\shf$, and 
$H^1\dsol(\DX^P,\shf)=\shf/P \shf$ describes the obstruction to solve 
the inhomogeneous equation $P\varphi=\psi$.

$\bullet$ Combining the two examples above, let 
$\nabla\colon\shg\to\shh$ be a differential operator acting between 
two holomorphic vector bundles.  Using the identification 
$\operatorname{Homdiff}(\shg,\shh) \simeq 
\Hom[\DX](\D\shh^*,\D\shg^*)$, we associate to $\nabla$ the 
$\DX$-module $\DX^\nabla$ represented by
$$
\D\shh^* \to \D\shg^* \to \DX^\nabla \to 0.
$$
One has $H^0\dsol(\DX^\nabla,\OX)=\ker\nabla$.

$\bullet$ Recall that a $\DX$-module is called holonomic if its 
characteristic variety has the smallest possible dimension, i.e.~is a 
Lagrangian subvariety of $T^*X$.  In other words, holonomic modules 
are locally associated to maximally overdetermined systems of PDE. If 
$\shm$ is holonomic, then $\dsol(\shm,\OX)$ is $\C$-constructible.  On 
the other hand, if $F$ is $\C$-constructible, then $\thom(F,\OX)$ has 
regular holonomic cohomology groups.  (Regularity is a generalization 
to higher dimensions of the notion of Fuchsian ordinary differential 
operators.)  More precisely, the Riemann-Hilbert correspondence 
asserts that the functors $\dsol(\cdot,\OX)$ and $\thom(\cdot,\OX)$ 
are quasi-inverse to each other.
\end{examples}

Recall that to a $d$-codimensional submanifold $S\subset X$ one
associates the regular holonomic $\D$-module of holomorphic
hyperfunctions $\B_{S|X}=H^d\thom(\C_{S},\OX)$.  Using more classical
notations, $\B_{S|X}= H^d_{[S]}(\OX)$ is a cohomology group with
tempered support.  The characteristic variety of $\B_{S|X}$ is the
conormal bundle $T^*_{S}X$.  If $S$ is a hypersurface, then elements
of $\B_{S|X}$ are equivalence classes of meromorphic functions with
poles along $S$, modulo holomorphic functions.  Even more
specifically, $\B_{\{0\}|\C}$ is generated by the equivalence class of
$(2\pi i z)^{-1}$, where $z\in\C$ is the holomorphic coordinate.

\subsection{Adjunction formula}
\label{ss:Adjunction}

In view of the above discussion, we consider the following setting for 
\eqref{eq:mfng}
\begin{align*}
    & \bigstar\text{ $X$ is a complex manifold} & & \bigstar\text{ 
    $\Xi$ is a complex manifold} \\
    & \bigstar\ \shm\in\BDC_{\coh}(\DX) & & \bigstar\ 
    \shn\in\BDC_{\coh}(\DY) \\
    & \bigstar\ F\in\BDC_{\rcons}(\C_{X}) & & \bigstar\ 
    G\in\BDC_{\rcons}(\C_{\Xi}) \\
    & \bigstar\text{ $\shf=\shc^\natural(F)$ for 
    $\natural=\pm\omega,\pm\infty$} & & \bigstar\text{ 
    $\shg=\shc^\natural(G)$ for $\natural=\pm\omega,\pm\infty$}
\end{align*}
Consider the projections
$$
X\from[q_{1}]X\times \Xi\to[q_{2}]\Xi.
$$
The three operations of pull-back, product, and push-forward, which 
appear in an expression like $\int\varphi(x)\cdot 
k(x,\xi)=\int_{q_{2}}(q_{1}^*\varphi\cdot k)$, make sense in the 
categories of sheaves and $\D$-modules.  We denote them by
$$
\opb{q_{1}},\ \tens,\ \reim{q_{2}}, \qquad\text{and}\qquad 
\dopb{q_{1}},\ \dtens,\ \deim{q_{2}},
$$
respectively.  For $\D$-modules, the analog of 
$\int_{q_{2}}(q_{1}^*\varphi\cdot k)$ is
$$
\shm\dcirc\shk = \deim{q_{2}}(\dopb{q_{1}}\shm\dtens\shk),
$$
where $\shm$ and $\shk$ are (complexes of) $\D$-modules on $X$ and on 
$X\times \Xi$, respectively.  Similarly, to the (complexes of) sheaves 
$G$ on $\Xi$ and $K$ on $X\times \Xi$, we associate the sheaf on $X$
$$
K\circ G = \reim{q_{1}}(\opb{q_{2}}G\tens K).
$$
Assume that $K$ and $\shk$ are interchanged by the Riemann-Hilbert 
correspondence, i.e.~assume one of the equivalent conditions
\begin{itemize}
    \item [$\bigstar$] $K$ is $\C$-constructible and 
    $\shk=\thom(K,\OXY)$,

    \item [$\bigstar$] $\shk$ is regular holonomic and 
    $K=\dsol(\shk,\OXY)$.
\end{itemize}
In order to ensure that $\R$-constructibility is preserved by the 
functor $K\circ \cdot$, and coherence is preserved by the functor 
$\cdot\dcirc\shk$, let us also assume that
\begin{itemize}
    \item [$\bigstar$] $q_1$ and $q_2$ are proper on $\supp(\shk)$,

    \item [$\bigstar$] $\shk$ is transversal to $\dopb{q_{1}}\shm$, 
    for any $\shm$,
\end{itemize}
where the latter condition means that $\chv(\shk)\inter(T^*X\times 
T^*_{\Xi}\Xi)$ is contained in the zero-section of $T^*(X\times \Xi)$.

Denote as usual by $(\cdot)\shift{1}$ the shift functor in derived 
categories.

\begin{theorem}
    \label{th:adj}
    With the above notations and hypotheses $\bigstar$, we have the 
    adjunction formula
    \begin{equation}
	\label{eq:adj}
	\DSol(\shm,\shc^\natural(K\circ G))\shift{\dim X} \simeq 
	\DSol(\shm\dcirc\shk,\shc^\natural(G)).
    \end{equation}
    In particular, to the pair of a $\D$-linear morphism 
    $\alpha\colon\shn\to\shm\dcirc\shk$ and a $\C$-linear morphism 
    $\beta\colon F\to K\circ G$ is associated a morphism $$
    \DSol(\shm,\shc^\natural(F))\shift{\dim X} \to[\beta\cdot\alpha] 
    \DSol(\shn,\shc^\natural(G)).  $$
\end{theorem}

This is the archetypical theorem we referred to in the Abstract.  In 
order to deal with concrete examples, one has to address three 
problems, which are of an independent nature:
\begin{itemize}
    \item [(T)] compute the sheaf transform $K\circ G$,

    \item [(A)] compute the $\D$-module transform 
    $\smash{\shm\dcirc\shk}$,

    \item [(Q)] identify the morphism $\beta\cdot\alpha$.
\end{itemize}

Problem (T) is purely topological, problem (A) is of an analytic 
nature, and (Q) is a quantization problem.  Usually the most difficult 
of the three is problem (A).  In the classical literature it 
corresponds to finding the system of PDE that characterizes the image 
of the transform.  This problem is split into two parts, of which the 
second is much more difficult
\begin{itemize}
    \item [(i)] find differential equations satisfied by functions in 
    the image of the transform,

    \item [(ii)] prove that these equations are sufficient to 
    characterize the image.
\end{itemize}
In our language, this translates into
\begin{itemize}
    \item [(i)] find a $\D$-linear morphism $\shn\to\shm\dcirc\shk$,

    \item [(ii)] prove it is an isomorphism.
\end{itemize}
An answer to the first issue is given by the following lemma.  
Although its proof is quite straightforward, it plays a key role.

\begin{lemma}
    \label{le:kernel}
    With the above notations and hypotheses $\bigstar$, one has an 
    isomorphism $$
    \alpha\colon H^0\DSol(\ddual\shm\detens\shn,\shk) \isoto 
    \Hom[\DY](\shn,\shm\dcirc\shk), $$
    where $\ddual\shm=\hom[\DX](\shm\tens[\OX]\Omega_{X},\DX)$ is the
    dual of $\shm$ as left $\DX$-modules, $\Omega_{X}$ denotes the
    sheaf of holomorphic forms of maximal degree, and
    $\Hom[\DY](\cdot,\cdot)$ is the set of globally defined
    $\DY$-linear morphisms.
\end{lemma}

This shows that any morphism $\shn\to\shm\dcirc\shk$ is of the form 
$\alpha(\kappa)$, where $\kappa$ is a globally defined $\shk$-valued 
solution of the system $\ddual\shm$ in the $X$ variables, and of the 
system $\shn$ in the $\Xi$ variables.  In other words, to any morphism 
$\shn\to\shm\dcirc\shk$ is attached an integral kernel $\kappa$.  We 
will show in Proposition~\ref{pr:quant} that, under some microlocal 
assumptions, one may as well read off from $\kappa$ the condition for 
$\alpha(\kappa)$ to be an isomorphism.

Concerning problem (Q), note that the morphism $\beta\cdot\alpha$ is 
obtained as the composite
\begin{eqnarray*}
    \DSol(\shm,\shc^\natural(F))\shift{\dim X} &\to& 
    \DSol(\shm,\shc^\natural(K\circ G))\shift{\dim X} \\
     &\simeq& \DSol(\shm\dcirc\shk,\shc^\natural(G))\\
     &\to& \DSol(\shn,\shc^\natural(G)),
\end{eqnarray*}
where the first and last morphisms are naturally induced by $\beta$ 
and $\alpha$, respectively.  One then has to solve the quantization 
problem of describing the distribution kernel $k$ of the integral 
transform $\beta\cdot\alpha$.  (Such kernels may be considered as the 
analog of Lagrangian distributions for Fourier Integral Operators.)  
In view of Lemma~\ref{le:kernel}, one should obtain $k$ as a boundary 
value of $\kappa$, $\beta$ describing the boundary value operation.  
An example of this phenomenon appears in Section~\ref{ss:realproj}, 
where we deal with the real projective Radon transform.

\subsection{Integral geometry}
\label{ss:IG}

In integral geometry one is given the graph $S\subset X\times \Xi$ of 
a correspondence from $\Xi$ to $X$.  Considering the natural morphisms
\begin{equation}
    X\from[f] S \to[g] \Xi
    \label{eq:XSY}
\end{equation}
induced by $q_{1}$ and $q_{2}$, the correspondence $\Xi \owns \xi 
\multimap \widehat \xi \subset X$ is given by $\widehat \xi=f(\opb g(\xi))$.  
For simplicity, let us assume that
\begin{itemize}
    \item [$\bigstar$] $S$ is a complex submanifold of $X\times \Xi$,

    \item [$\bigstar$] $f$ and $g$ are smooth and proper.
\end{itemize}
One is interested in those transforms whose integral kernel $k$ is 
some characteristic class of $S$, so that 
$\int_{q_{2}}(q_{1}^*\varphi\cdot k)=\int_{\widehat \xi}\varphi$ becomes 
an integration along the family of subsets $\smash{\widehat \xi}$.  More 
precisely, Theorem~\ref{th:adj} fits in the framework of integral 
geometry if $\shk$ is a regular holonomic module whose characteristic 
variety coincides---at least outside of the zero section---with the 
conormal bundle $T^*_{S}(X\times \Xi)$.

In this context, a natural choice is $\shk=\B_{S}$, where we use the 
shorthand notation $\B_{S}=\B_{S|X\times \Xi}$ for the sheaf of 
holomorphic hyperfunctions.  The Riemann-Hilbert correspondence 
associates to $\shk=\B_{S}$ the perverse complex $K=\C_{S}\shift{-d}$, 
where $d=\codim_{X\times \Xi}S$.  One easily checks that
$$
K\circ G\simeq\reim f\opb g G\shift{-d}, \qquad 
\shm\dcirc\shk\simeq\deim g\dopb f\shm.
$$

Let us show in which sense problem (T) of the previous section is 
purely topological.  For example, let $G=\C_{Z}$ for $Z\subset \Xi$ a 
locally closed subanalytic subset.  Then, the fiber of $K\circ G$ at 
$x\in X$ is computed by
\begin{align}
    \label{eq:fiber}
    (\reim f\opb g \C_{Z})_{x} &\simeq (\reim f \C_{\opb g(Z)})_{x}\\
    \nonumber &\simeq \rsect_{c}(\opb f(x); \C_{\opb f(x) \inter \opb 
    g(Z)})\\
    \nonumber &\simeq \rsect_{c}(\widehat x; \C_{\widehat x \inter Z}),
\end{align}
where in the last isomorphism we used the identification $g\colon\opb 
f(x) \inter \opb g(Z)\isoto\widehat x \inter Z$.  This shows that $K\circ 
G$ is determined by the (co)homology groups of the family of slices 
$x\multimap\widehat x \inter Z$.

\subsection{Microlocal geometry}
\label{ss:Microlocal}

Let us consider the microlocal correspondence associated to 
\eqref{eq:XSY}
\begin{equation}
    T^*X \from[p_{1}] T^*_{S}(X\times \Xi) \to[p_{2}^a] T^*\Xi,
    \label{eq:muXSY}
\end{equation}
where $p_{1}$ and $p_{2}$ are induced by the natural projections from 
$T^*(X\times \Xi)$, and $(\cdot)^a$ denotes the antipodal map of 
$T^*\Xi$.  By restriction, this gives the correspondence
\begin{equation}
    \dot T^*X \from[\dot p_{1}] \dot T^*_{S}(X\times \Xi) \to[\dot 
    p_{2}^a] \dot T^*\Xi,
    \label{eq:dotmuXSY}
\end{equation}
where the dot means that we have removed the zero-sections.

Using the theory of microdifferential operators, we may get some {\em 
a priori} information on $\shm\dcirc\B_{S}$.

\begin{proposition}
    \label{pr:est}
    With the above notations and hypotheses $\bigstar$, assume that 
    the map $\dot p_{2}^a$ is finite.  Then $$
    \chv(\shm\dcirc\B_{S})\subset p_{2}^a(\opb{p_{1}}(\chv(\shm))).  
    $$
    Moreover, if $\shm$ is concentrated in degree zero, then 
    $H^j(\shm\dcirc\B_{S})$ is a flat connection for $j\neq 0$.  
    Finally, if $\shh$ is a holomorphic vector bundle, then $H^j(\D 
    \shh^* \dcirc \B_{S}) = 0$ for $j<0$.
\end{proposition}

Assume now
\begin{itemize}
    \item [$\bigstar$] $\dot p_{1}$ is smooth surjective,

    \item [$\bigstar$] $\dot p_{2}^a$ is a closed embedding onto a 
    regular involutive submanifold $\dot \MV\subset\dot T^*\Xi$.
\end{itemize}
The above conditions imply that the fibers of $p_{1}$ are identified 
with the bicharacteristic leaves of $\dot \MV$.  Thus, it is not 
difficult to prove that, locally on $\dot T^*_{S}(X\times \Xi)$, the 
correspondence \eqref{eq:dotmuXSY} is isomorphic to a contact 
transformation with parameters.

Let us recall three notions from the theory of microdifferential
operators.  First, one says that a $\DY$-module $\shn$ has simple
characteristics along $\dot\MV$ if in a neighborhood of any point in
$\dot\MV$, its associated microdifferential system admits a simple
generator whose symbol ideal is reduced and coincides with the
annihilating ideal of $\dot \MV$.  In particular, if $\shh$ is a line
bundle on $X$, then $\D\shh^*$ is simple along $\dot T^*X$.  Second,
there is a natural notion of non-degenerate section of $\B_{S}$, which
is local on $\dot T^*_S(X\times \Xi)$.  For example, if $S$ is a
smooth hypersurface with local equation $(\Phi=0)$, then a section of
$\B_{S}$ is non-degenerate if it is obtained by applying an invertible
microdifferential operator to the generator $1/\Phi$.  Finally, one
says that a $\D$-linear morphism is an m-f-c isomorphism (short hand
for isomorphism modulo flat connections) if its kernel and cokernel
are locally free $\O$-modules of finite rank.

\begin{proposition}
    \label{pr:quant}
    With the above notations and hypotheses $\bigstar$, let $\shh$ be 
    a holomorphic vector bundle on $X$, and let $\shn$ be a 
    $\DY$-module with simple characteristics along $\dot\MV$.  Let 
    $\kappa\in H^0\DSol(\ddual{(\D \shh^*)}\detens\shn,\shk)$ be a 
    non-degenerate section of $\B_{S}$.  Then, 
    $$H^0\alpha(\kappa)\colon \shn \to H^0(\D \shh^*\dcirc\B_{S})$$ is 
    an m-f-c isomorphism.
\end{proposition}

Applying \eqref{eq:adj} with $G=\C_{\xi}$, $\xi\in \Xi$, and 
$\natural=\omega$, we get the germ formula
\begin{equation}
    \label{eq:germ}
    (\D\shh^*\dcirc\B_{S})_{\xi} \simeq
    \rsect(\widehat\xi;\shh\vert_{\widehat\xi})\shift{\dim S-\dim
    \Xi}.
\end{equation}
This, together with the similar formula for $\ddual{(\D\shh^*)}$,
provide a useful test to check both if $\alpha(\kappa) =
H^0\alpha(\kappa)$ and if $\alpha(\kappa)$ is an actual isomorphism,
not only modulo flat connections.

One says that a $\DY$-module $\shn$ has regular singularities along 
$\dot \MV$ if it locally admits a presentation
$$
\shs^{N_{1}}\to \shs^{N_{0}} \to \shn \to 0,
$$
where $\shs$ has simple characteristics along $\dot\MV$.  Assume
\begin{itemize}
    \item [$\bigstar$] $f$ has connected and simply connected fibers.
\end{itemize}

\begin{theorem}
    \label{th:equivalence}
    With the above notations and hypotheses $\bigstar$, the functor 
    $\shm\mapsto H^0(\shm\dcirc\B_{S})$ induces an equivalence of 
    categories between coherent $\DX$-modules, modulo flat 
    connections, and $\DY$-modules with regular singularities along 
    $\dot \MV$, modulo flat connections.  Moreover, this equivalence 
    interchanges modules with simple characteristics along $\dot T^*X$ 
    with modules with simple characteristics along $\dot \MV$.  
    Finally, if $\dim X\geq 3$, any $\DX$-module with simple 
    characteristics along $\dot T^*X$ is m-f-c isomorphic to a module 
    of the form $\D\shh^*$, for some line bundle $\shh$ on $X$.
\end{theorem}

\subsection*{Notes}

{\bf \S\ref{ss:GeneralizedFcts}} The idea of constructing generalized
functions starting with holomorphic functions is at the heart of
Sato's theory of hyperfunctions.  When making operations on
hyperfunctions one is then led to consider complexes of the form
$\shc^{\pm\omega}(F)$.  This theory is developed
in~\cite{Kashiwara-Schapira90}, which is also a very good reference on
sheaf theory in derived category.

\noindent By taking into account growth conditions, distributions can 
also be obtained from the sheaf of holomorphic functions.  This is 
done using the functor $\thom$ of tempered cohomology introduced 
in~\cite{Kashiwara84}.  The dual construction of the formal cohomology 
functor $\wtens$ is performed in~\cite{Kashiwara-Schapira96}, where 
one also finds a systematic treatment of operations on 
$\dsol(\shm,\shc^{\pm\infty}(F))$.

\medskip\noindent
{\bf \S\ref{ss:PDE}} Kashiwara's Master Thesis, recently translated 
in~\cite{Kashiwara95}, is still a very good reference on the analytic 
theory of $\D$-modules.  See also~\cite{Sato-Kawai-Kashiwara73}, 
\cite{Bjork93} and~\cite{Schneiders94}.  A proof of the 
Riemann-Hilbert correspondence is obtained in~\cite{Kashiwara84} by 
showing that the functor $\thom(\cdot,\OX)$ is a quasi-inverse to the 
solution functor $\dsol(\cdot,\OX)$.

\medskip\noindent {\bf \S\ref{ss:Adjunction}} The language of
correspondences is very classical.  In the category of sheaves, it can
be found for example in~\cite{Kashiwara-Schapira90}.  The formalism of
$\dsol(\shm,\shc^{\pm\omega}(F))$ was introduced
in~\cite{Schapira-Schneiders94}.  The idea of using this framework to
investigate integral transforms is from~\cite{D'Agnolo-Schapira96a}.

\noindent The starting point to address problems in integral geometry
is the adjunction formula in Theorem~\ref{th:adj}.  For
$\natural=\pm\omega$ this was obtained in~\cite{D'Agnolo-Schapira96a},
\cite{D'Agnolo-Schapira96b} using results as the
Cauchy-Kovalevskaya-Kashiwara theorem and Schneiders' relative duality
theorem.  The case $\natural=\pm\infty$ is
from~\cite{Kashiwara-Schapira96}.  Note also
that~\cite{Kashiwara-Schmid94} independently announced a similar
result in an equivariant framework, at the time of our
announcement~[A.~D'Agnolo and P.~Schapira, C.~R.\ Acad.\ Sci.\ Paris
S{\'e}r.\ I Math.\ {\bf 319}, no.~5 (1994), 461--466; Ibid.\ no.~6,
595--598].

\noindent The kernel lemma~\ref{le:kernel} is 
from~\cite{D'Agnolo-Schapira96b} (see also~\cite{Goncharov97} 
and~\cite{Tanisaki98}).

\noindent We refer to the appendix of~\cite{D'Agnolo98} for a 
discussion on how to obtain the distribution kernel $k$ as boundary 
value of the meromorphic kernel $\kappa$.

\medskip\noindent
{\bf \S\ref{ss:Microlocal}} The microlocal geometry attached to double 
fibrations was first considered in~\cite{Guillemin-Sternberg79}.
Proposition~\ref{pr:est} is 
from~\cite{D'Agnolo-Schapira96a}.  Its proof relies on a result 
asserting that the functor $\cdot\dcirc\B_{S}\simeq\deim g\dopb 
f(\cdot)$ commutes to microlocalization.  The corresponding result for 
inverse images is due to~\cite{Sato-Kawai-Kashiwara73}, and the one 
for direct images is due to~\cite{Schapira-Schneiders94}.

\noindent Proposition~\ref{pr:quant} is 
from~\cite{D'Agnolo-Schapira96b}, \cite{D'Agnolo-Schapira98}.  Its 
proof is based on the quantization of contact transformations 
by~\cite{Sato-Kawai-Kashiwara73}, where the notion of non degenerate 
section is found (see also~\cite{Schapira85} for an exposition of the 
theory of microdifferential operators).
The equivalence result in Theorem~\ref{th:equivalence} was 
obtained in~\cite{D'Agnolo-Schapira96a}, \cite{D'Agnolo-Schapira98}.  
The case $X=\PP$, $\Xi=\bb$ had already been considered by 
Brylinski~\cite{Brylinski86} for regular holonomlic modules.

\section{Back to the Radon transform}
\label{se:proofs}

Using the framework of sheaves and $\D$-modules, we will give here a 
proof of the statements in Section~\ref{se:Statements}, along with 
some generalizations.

According to \eqref{eq:geoC}, we denote by $\PP$ the projective space 
of vector lines in $\CV\simeq\C^{n+1}$, and by $\GG$ the Grassmannian 
of $(p+1)$-dimensional subspaces of $\CV$.  Let 
$\FF\subset\PP\times\GG$ be the incidence relation of pairs 
$(z,\zeta)$ with $z\in\widehat\zeta$.  Recall that $\dim\PP=n$, 
$\dim\GG=(p+1)(n-p)$, and $\FF$ is a smooth submanifold of 
$\PP\times\GG$ of dimension $n+p(n-p)$.  By definition, $\FF$ is the 
graph
\begin{equation}
    \label{eq:PFG}
    \PP \from[f] \FF \to[g] \GG
\end{equation}
of the correspondence $\GG\owns\zeta\multimap\widehat\zeta\subset\PP$ 
attached to the complex Radon transform.  Let us describe the 
associated microlocal correspondence
\begin{equation}
    \label{eq:muPFG}
    T^*\PP \from[p_{1}] T^*_{\FF}(\PP\times\GG) \to[p_{2}^a] T^*\GG.
\end{equation}
Denote by $\langle\zeta\rangle\subset \CV$ the $(p+1)$-vector subspace 
attached to $\zeta\in\GG$.  It is a nice exercise in classical 
geometry\footnote{see for example [J.\ Harris, {\sl Algebraic 
geometry.  A first course.} Graduate Texts in Mathematics, 133.  
Springer-Verlag, New York, 1995, \MR{97e:14001}]} to recover the 
formula $(T\GG)_{\zeta} = 
\Hom(\langle\zeta\rangle,\CV/\langle\zeta\rangle)$, where this $\Hom$ 
is simply the set of linear morphisms of vector spaces.  Similarly, 
one has
\begin{eqnarray*}
     T^*\PP & = & \{ (z;\alpha)\colon \alpha \in \Hom(\CV/\langle 
     z\rangle,\langle z\rangle)\}, \\
     T^*\GG & = & \{ (\zeta;\beta)\colon \beta \in 
     \Hom(\CV/\langle\zeta\rangle,\langle\zeta\rangle)\}, \\
     T^*_{\FF}(\PP\times\GG) & = & \{ (z,\zeta;\gamma)\colon \gamma 
     \in \Hom(\CV/\langle\zeta\rangle,\langle z\rangle)\}.
\end{eqnarray*}
The maps \eqref{eq:muPFG} are then given by
$$
p_{1}(z,\zeta;\gamma)=(z;\gamma\circ q), \qquad 
p_{2}^a(z,\zeta;\gamma)=(\zeta;j\circ \gamma),
$$
where $q\colon \CV/\langle z\rangle \twoheadrightarrow 
\CV/\langle\zeta\rangle$ and $j\colon \langle z\rangle \hookrightarrow 
\langle\zeta\rangle$ are induced by the inclusion $\langle z\rangle 
\subset \langle\zeta\rangle$.

\begin{remark}
    One should be careful that, with the notations \eqref{eq:zeta}, $$
    \langle\zeta\rangle = 
    (\C\cdot\zeta_{1}+\cdots+\C\cdot\zeta_{n-p})^\bot \subset \CV. $$
    This is due to the fact that \eqref{eq:zeta} uses the 
    identification $\zeta\mapsto\zeta^\bot$ of $\GG$ with the 
    Grassmannian of $(n-p)$-dimensional subspaces of $\CV^*$.
\end{remark}

It is now easy to check that all of the hypotheses $\bigstar$ in the 
previous section are satisfied for the choice
\begin{itemize}
    \item [$\blacktriangle$] $X=\PP$,\quad $\Xi=\GG$,\quad 
    $S=\FF$,\quad $\shk=\BF$,\quad $K=\C_{\FF}\shift{p-n}$.
\end{itemize}
Note also that $\MV=p_{2}^a(p_{1}^{-1}(T^*\PP))$ is described by
\begin{itemize}
    \item [$\blacktriangle$] $\MV=\{(\zeta;\beta)\in T^*\GG\colon
    \beta\colon\CV/\langle\zeta\rangle\to\langle\zeta\rangle \text{
    has rank at most one}\}$.
\end{itemize}
By Proposition~\ref{pr:est}~(i), $\MV$ is the {\em a priori} estimate
for the characteristic variety of any $\DG$-module that can possibly
arise from the Radon transform.  Not too surprisingly, this is the
characteristic variety of the Maxwell-John system $\square$.

Our aim here is to show that all of the examples from 
Section~\ref{se:Statements} are particular cases of 
Theorem~\ref{th:adj}, for different choices of $\alpha \colon \shn \to 
\shm\dcirc\BF$ and $\beta \colon F \to \C_{\FF}\circ G\shift{p-n}$.  
Even better, for those examples $\alpha$ is always the same.

\subsection{A Radon $\D$-module}
\label{se:RadonD}

The real affine Radon transform in Theorem~\ref{th:Rreal} deals with 
the space $\shs(\RA)$, where no differential equations appear.  
Theorems~\ref{th:Rcomplex} and~\ref{th:hyp} are concerned with 
function spaces attached to the line bundle $\OP\twist{-p-1}$.  Again, 
this is something which is locally trivial.  It is then natural to 
choose $\shm=\DP\twist{p+1}$, where we set 
$\DP\twist{m}=\DP\tens[\OP]\OP\twist{m}$ for $m\in\Z$.

Concerning $\shn$, if $p<n-1$ this should be the $\DG$-module 
$\DG^{\square}$ represented by the Maxwell-John system \eqref{eq:wave}
$$
\D\shh^*\to[\square]\DG\twist{1}\to\DG^{\square}\to 0.
$$
If $p=n-1$, then $\GG=\bb$ is a dual projective space and we should
take $\shn=\Db\twist{1}$.  Since in this case $\square=0$, we may
still write $\Db\twist{1} = \Db^{\square}$.

Finally, the choice of $\kappa$ is imposed by \eqref{eq:kRC}.  
Summarizing, we consider
\begin{itemize}
    \item [$\blacktriangle$] $\shm=\DP\twist{p+1}$,\quad 
    $\shn=\DG^{\square}$,\quad $(2\pi i)^{n-p}\kappa(z,\zeta) = 
    \dfrac{\omega(z)}{\langle z,\zeta_{1}\rangle \cdots \langle 
    z,\zeta_{n-p}\rangle}$.
\end{itemize}
To explain the meaning of $\kappa$, note that the \v Cech covering 
$\{\langle z,\zeta_{i}\rangle\neq 0\}_{i=1,\dots,n-p}$ of 
$(\PP\times\GG) \setminus \FF$ allows us to locally consider $\kappa$ 
as a cohomology class in $H^{n-p}_{[\FF]}\OPG$.  Globally, 
$\kappa(z,\zeta)$ is $(p+1)$-homogeneous in $z$, and 
$(-1)$-homogeneous in $\zeta$, for the action of $\GL(1,\C)$ and 
$\GL(n-p,\C)$, respectively.  This is written as
$$
\kappa\in H^0\DSol(\ddual{\DP\twist{p+1}}\detens\DG\twist{1},\BF).
$$
Moreover, it is easy to check that $\kappa(z,\zeta)$ is a solution of 
$\square$, acting on the $\zeta$ variable.  Hence
$$
\kappa\in H^0\DSol(\ddual{\DP\twist{p+1}}\detens\DG^{\square},\BF).
$$
It is now a local problem on $\dot T^*_{\FF}(\PP\times\GG)\isoto\dot 
\MV$ to verify that $\DG^{\square}$ is simple along $\dot \MV$, and 
that $\kappa$ is non-degenerate.

\begin{theorem}
    \label{th:kappa-p-1}
    The above choice of $\kappa$ induces an isomorphism $$
    \alpha(\kappa)\colon\DG^{\square}\isoto\DP\twist{p+1}\dcirc\BF. $$
\end{theorem}

\begin{remark}
    The above discussion was facilitated by the fact that the choice
    of $\shm$, $\shn$ and $\kappa$ was forced by the statements we
    wanted to recover.  In general, one either has good candidates for
    $\shn$ and $\shk$, and may then proceed as above, or one has to
    actually compute the transform $\shm\dcirc\shk$.  This last
    problem is generally difficult.  Note, however, that many of the
    examples that arise in practice are endowed with the action of a
    group, including the Radon transform which is equivariant for the
    action of $\GL(\CV)$.  Taking this into account can greatly
    simplifies matters, by narrowing down the possible outcome of the
    transform, and by allowing one to use the computational techniques
    of representation theory.
\end{remark}

\subsection{Radon adjunction formula}

For an $\OP$-module $\shf$, set 
$\shf\twist{m}=\shf\tens[\OP]\OP\twist{m}$.  Combining 
Theorems~\ref{th:adj} and~\ref{th:kappa-p-1} we get

\begin{theorem}
    \label{th:Radon}
    Let $\beta\colon F\to \C_{\FF}\circ G$ be a $\C$-linear morphism 
    inducing an isomorphism
    \begin{equation}
	\label{eq:piso}
	H^p(\PP;\shc^\natural(F)\twist{-p-1}) \isoto 
	H^p(\PP;\shc^\natural(\C_{\FF}\circ G)\twist{-p-1}).
    \end{equation}
    Then, one has an isomorphism $$
    \beta\cdot\alpha(\kappa)\colon 
    H^p(\PP;\shc^\natural(F)\twist{-p-1}) \isoto 
    H^0\DSol(\DG^{\square},\shc^\natural(G)).  $$
\end{theorem}

Requiring that $\beta$ itself be an isomorphism is of course a
sufficient condition for \eqref{eq:piso} to hold, but it is not
necessary.  For example, since $\rsect(\PP;\OP\twist{-p-1})= 0$, one
has that $\rsect(\PP;\shc^\natural(N)\twist{-p-1})=0$ if $N$ is a
complex of finite rank constant sheaves on $\PP$.  In other words, the
morphism $\rsect(\PP;\shc^\natural(\beta)\twist{-p-1})$ is an
isomorphism if $\beta$ is an isomorphism in the localization of
$\BDC_{\rcons}(\C_{\PP})$ by the null systems of objects like $N$.  We
denote this localization by $\BDC_{\rcons}(\C_{\PP};\dot T^*\PP)$.

We are now in a position to deduce all of the results in 
Section~\ref{se:Statements} from Theorem~\ref{th:Radon}, for different 
choices of $\beta\colon F\to \C_{\FF}\circ G$.

\subsection{Complex projective case}
\label{ss:complexproj2}

We use the same notations as in Section~\ref{ss:complexproj}.  Recall 
that
$$
\shc^{-\omega}(\C_{\widehat U}')=\rsect_{\widehat U}(\OP), \quad 
\shc^{-\omega}(\C_{U}')=\rsect_{U}(\OG).
$$
To get Theorem~\ref{th:Rcomplex} as a corollary of 
Theorem~\ref{th:Radon}, one should then take
\begin{itemize}
    \item [$\blacktriangle$] $F=\C_{\widehat U}'$,\quad $G=\C_{U}'$,\quad 
    $\natural=-\omega$.
\end{itemize}
Note that, by duality, the datum of a morphism
$$
\beta\colon F\to \C_{\FF}\circ G
$$
is equivalent to the datum of
$$
\beta'\colon F'\from (\C_{\FF}\circ G)' \simeq \C_{\FF}\circ 
G'\shift{2p(n-p)}.
$$
According to~\eqref{eq:fiber}, for $z\in\PP$ one has
\begin{align*}
    H^j(\C_{\FF}\circ\C_{U})_{z} &\simeq H^j_{c}(\widehat z;\C_{\widehat z 
    \inter U}) \\
    &\simeq H^{2p(n-p)-j}(\widehat z;\C_{\widehat z \inter U}),
\end{align*}
where the last isomorphism is Poincar{\'e} duality.  The set $\widehat z 
\inter U$ is non-empty if and only if $z\in\widehat U$.  If $\widehat z \inter 
U$ is connected for any $z\in\widehat U$, we get a morphism
$$
\beta'\colon \C_{\widehat U}\from \C_{\FF}\circ \C_{U}\shift{2p(n-p)}.
$$
Since $U$ is elementary, the truncation $\tau^{\geq -p}(\beta')$ of 
$\beta'$ in degree greater or equal to $-p$ is an isomorphism.  From 
this fact one deduces \eqref{eq:piso}, and the statement follows.

\begin{remark}
    The technical point of passing from $\beta$ to $\beta'$ is 
    inessential and could have been avoided.  It is solely due to our 
    definition of $\shc^{-\omega}(F)$, which incorporates $F'$ to be 
    best suited for the real case.
\end{remark}

Of course, we could also consider the cases 
$\natural=\omega,\pm\infty$.  For example, take
\begin{itemize}
    \item [$\blacktriangle$] $p=n-1$,\quad $F=\C_{\widehat U}$,\quad 
    $G=\C_{U}\shift{2(n-1)}$,\quad $\natural=\omega$.
\end{itemize}
If $U$ is a bounded neighborhood of the origin in an affine chart 
$\aa\subset\bb$, then $U^\sharp=\PP\setminus\widehat U$ is a compact 
subset in an affine chart $\AA\subset\PP$.  We thus recover 
Martineau's isomorphism
$$
\sect(U^\sharp;\O_\AA) \simeq \O'_{\aa}(U),
$$
where $\O'_{\aa}(U) = H^n_{c}(U;\O_{\aa}) = H^n_{c}(U;\Ob\twist{-1})$ 
is the space of analytic functionals in $U$.

\subsection{Real conformal case}
\label{se:realconf2}

We use the same notations as in Section~\ref{ss:realconformal}.  By 
definition, we have
$$
\rsect(\PP;\shc^\omega(\C_{\UP})) \simeq \rsect(\UP;\OP), \quad 
\rsect(\PP;\shc^{-\omega}(\C_{\UP})) \simeq 
\rsect_{\UP}(\PP;\OP)\shift{1},
$$
where in the second isomorphism we used the identification 
$\C_{\UP}'\simeq\C_{\UP}\shift{1}$, due to the fact that $\UP$ is a 
smooth hypersurface splitting $\PP$ in two connected components.  
Since $\UG$ is totally real in $\GG$, we also have
$$
\shc^\omega(\C_{\UG}) = \shc^\omega_{\UG}, \quad 
\shc^{-\omega}(\C_{\UG}) = \shc^{-\omega}_{\UG}.
$$
Take
\begin{itemize}
    \item [$\blacktriangle$] $F=\C_{\UP}$,\quad $G=\C_{\UG}$,\quad 
    $\natural=\pm\omega$,
\end{itemize}
and let $\beta$ be the morphism induced by the equality $\widehat\UG = 
\UP$.  Then, the choice $\natural=\omega$ allows to get the 
isomorphism in the first line of Theorem~\ref{th:hyp}, while 
$\natural=-\omega$ gives the one in the second line.

\subsection{Real projective case}
\label{ss:realproj}

Let us discuss this case in somewhat more detail.  Consider
\begin{equation}
    \begin{cases}
	\RV\simeq\R^{n+1} &\text{a real vector space with 
	$\RV\tens[\R]\C=\CV$,} \\
	\RP &\text{the projective space of vector lines in $\RV$,} \\
	\RG &\text{the Grassmannian of projective $p$-planes in 
	$\RP$.}
    \end{cases}
    \label{eq:geoRP}
\end{equation}
Since $\pi_{1}(\RP)=\Z/2\Z$ (let us assume here that $n>1$), there are 
essentially two locally constant sheaves of rank one on $\RP$.  For 
$\varepsilon\in\Z/2\Z$, we denote them by $\C_{\RP}^{(\varepsilon)}$, 
asking that $\C_{\RP}^{(0)}$ be the constant sheaf $\C_{\RP}$.  One 
then easily checks that
$$
\shc^\infty(\C_{\RP}^{(-p-1)})\twist{-p-1} \simeq 
\shc^{\infty}_{\RP}\abstwist{-p-1},
$$
where the term on the right hand side denotes the $C^\infty$ line 
bundle on $\RP$ whose sections $\varphi$ satisfy the homogeneity 
condition
$$
\varphi(\lambda x)=|\lambda|^{-p-1} \varphi(x) \quad 
\forall\lambda\in\GL(1,\R).
$$
Since $\pi_{1}(\RG)=\Z/2\Z$, in exactly the same way we have
$$
\quad \shc^{\infty}_{\RG}\abstwist{-1} = 
\shc^\infty(\C_{\RG}^{(-1)})\twist{-1}.
$$
(Note that if $q\colon\RG^+\to\RG$ denotes the $2:1$ projection from 
the Grassmannian of oriented planes to $\RG$, one has $\oim 
q(\C_{\RG^+})=\C_{\RG}\dsum\C_{\RG}^{(1)}$.)

Let us set
\begin{itemize}
    \item [$\blacktriangle$] $F=\C_{\RP}^{(-p-1)}$,\quad 
    $G=\C_{\RG}^{(-1)}$,\quad $\natural=\infty$.
\end{itemize}
Denote by $\RG_{p,n}$ the Grassmannian of $p$-dimensional subspaces of 
$\R^n$.  For $z\in\PP$ one has
$$
\widehat z\inter \RG = \{\xi\in\RG \colon \C\tens[\R]\langle\xi\rangle 
\supset \langle z \rangle\} \simeq
\begin{cases}
    \RG_{p,n}, & \text{for }z\in\RP, \\
    \RG_{p-1,n-1}, & \text{for }z\in\PP\setminus\RP.
\end{cases}
$$
By a computation like \eqref{eq:fiber}, this implies
$$
H^j(\C_{\FF}\circ\C_{\RG}^{(-1)})_{z} \simeq
\begin{cases}
    H^j(\RG_{p,n};\C_{\RG_{p,n}}^{(-1)}), & \text{for }z\in\RP, \\
    H^j(\RG_{p-1,n-1};\C_{\RG_{p-1,n-1}}^{(-1)}), & \text{for 
    }z\in\PP\setminus\RP.
\end{cases}
$$
Looking at a table of Betti numbers for oriented and non-oriented real 
Grassmannians, one deduces that
$$
\tau^{\leq p}(\C_{\FF}\circ\C_{\RG}^{(-1)}) \simeq
\begin{cases}
    \C_{\RP}^{(-p-1)}, & \text{for $p$ even}, \\
    \C_{\PP\setminus\RP}, & \text{for $p$ odd}.
\end{cases}
$$
Denoting by
$$
\beta\colon \C_{\RP}^{(-p-1)}\to\C_{\FF}\circ\C_{\RG}^{(-1)}
$$
the natural morphism, the above arguments imply that $\tau^{\leq 
p}(\beta)$ is an isomorphism in $\BDC_{\rcons}(\C_{\PP};\dot T^*\PP)$.  
Applying Theorem~\ref{th:Radon}, we get a compactified version of 
Theorem~\ref{th:Rreal}.

\begin{theorem}
    \label{th:Rproj}
    The real projective Radon transform
    \begin{align*}
	R_{\RP} \colon \sect(\RP;\shc^\infty_{\RP}\abstwist{-p-1}) 
	&\to \sect(\RG;\shc^\infty_{\RG}\abstwist{-1}) \\
	\varphi(x) &\mapsto \psi(\xi) = \int \varphi(x) \delta(\langle 
	x,\xi_{1}\rangle) \cdots \delta(\langle x,\xi_{n-p}\rangle) 
	\omega(x)
    \end{align*}
    induces an isomorphism $$
    R_{\RP} \colon \sect(\RP;\shc^\infty_{\RP}\abstwist{-p-1}) \isoto 
    \sect(\RG;\ker(\square,\shc^\infty_{\RG})).  $$
\end{theorem}
Here, we used the identification $R_{\RP}=\beta\cdot\alpha(\kappa)$, 
which follows from the fact that
$$
k(x,\xi) = \delta(\langle x,\xi_{1}\rangle) \cdots \delta(\langle 
x,\xi_{n-p}\rangle) \omega(x)
$$
is the boundary value of $\kappa(z,\zeta)$.

\subsection{Real affine case}
\label{ss:realaffine2}

The affine case \eqref{eq:geoRA} sits in the projective case 
\eqref{eq:geoRP} by considering
\begin{equation}
    \begin{cases}
	\RA=\RP\setminus\RH &\text{for a hyperplane $\RH\subset\RP$,} 
	\\
	\RGA\subset\RG &\text{the set of $\xi\in\RG$ with 
	$\widehat\xi\subset\RA$.}
    \end{cases}
    \label{eq:geoRAinRP}
\end{equation}
The space $\shs(\RA)$ is then identified to the space of 
$C^\infty$-functions globally defined in $\RP$, which vanish up to 
infinite order on $\RH$.  Since twisting does not matter in affine 
charts, we have
$$
\shs(\RA) \simeq \sect(\PP;\shc^\infty(\C_{\RA})\twist{-p-1}).
$$
To recover Theorem~\ref{th:Rreal}, one then has to consider
\begin{itemize}
    \item [$\blacktriangle$] $F=\C_{\RA}\simeq\C_{\RA}^{(-p-1)}$,\quad 
    $G=\C_{\RGA}^{(-1)}$,\quad $\natural=\infty$,
\end{itemize}
where $\C_{\RGA}^{(-1)}=\C_{\RGA}\tens\C_{\RG}^{(-1)}$.  Let us detail 
the more interesting case $p=n-1$, where the Cavalieri condition 
appears.

We are now considering the situation
\begin{equation}
    \begin{cases}
	\RG=\Rb &\text{a real projective space in $\GG=\bb$,} \\
	\RGA=\Rbo &\text{where $\Rbo=\Rb\setminus\{\xio\}$ for 
	$\xio\in\Rb$ with $\widehat\xio=\RH$.}
    \end{cases}
    \label{eq:geoRAinRPforRb}
\end{equation}
Since the Radon hyperplane transform is symmetrical, we may apply 
Theorem~\ref{th:Radon} interchanging the roles of $\PP$ and $\bb$.  We 
get
\begin{equation}
    \label{eq:imRa}
    R_{\RA}(\shs(\RA)) \simeq 
    H^n(\bb;\shc^\infty(\C_{\RA}\circ\C_{\FF})\twist{-1}).
\end{equation}

The embedding $\RH \subset \RP$ induces a projection $\Rbo\to\Rh$, and 
its complexification
$$
q\colon \bbo \to \hh,
$$
where $\Rh$ is the dual projective space to $\RH$.  Note that the 
fibers of $q$ are the complex projective lines through $\xio$, with 
the point $\xio$ removed.  Consider
\begin{equation}
    \label{eq:qRb}
    \qRbo = \opb q (\Rh), \quad \qRb = \qRbo \union \{\xio\}.
\end{equation}
The set $\qRb$ has its only singularity at $\xio$.  An explicit 
computation gives the distinguished triangle in 
$\BDC_{\rcons}(\C_{\PP};\dot T^*\PP)$
$$
\C_{\RA}\circ\C_{\FF} \to \C_{\Rbo}^{(-1)}\shift{-n} \to 
\C_{\qRb}\shift{1-n} \to[+1].
$$
Using the notation $\shs_{(-1)}(\qRb) = 
\rsect(\bb;\shc^\infty(\C_{\qRb})\twist{-1})$,
we may then write \eqref{eq:imRa} as
$$
R_{\RA}(\shs(\RA)) = \ker \left(\, \shs_{\{-1\}}(\Rbo) \to[c] 
H^1\shs_{(-1)}(\qRb) \, \right),
$$
and we are left to describe the morphism $c$.

The decomposition \eqref{eq:qRb} gives an identification
$$
H^1 \shs_{(-1)}(\qRb) = \coker \left(\, \Ob \twist{-1} 
\widehat\vert_{\xio} \to[j] H^1 \shs_{(-1)}(\qRbo) \, \right).
$$
Here $\Ob\twist{-1}\widehat\vert_{\xio} = H^1\shs_{(-1)}(\{\xio\})$ is 
the formal restriction of $\Ob\twist{-1}$ to $\xio$, whose sections 
are formal Taylor series
$$
\Ob\twist{-1}\widehat\vert_{\xio} \simeq \prod_{m\geq 
0}\sect(\hh;\O_{\hh}\twist{m}).
$$
Summarizing, we have
$$
\shs_{\{-1\}}(\Rbo) \to[\tilde c] \prod_{m\geq 
0}\sect(\Rh;\shc^\infty_{\Rh}\twist{m}) \isofrom 
H^1\shs_{(-1)}(\qRbo) \from[j]\prod_{m\geq 
0}\sect(\hh;\O_{\hh}\twist{m}),
$$
where $\tilde c$ and the middle isomorphism are obtained by 
integration along the fibers of $q$.  Using the above identifications, 
we finally get
$$
R_{\RA}(\shs(\RA)) = \{\psi \in \shs_{\{-1\}}(\Rbo) \colon \tilde 
c(\psi)_{m} = j(\Psi_{m}) \text{ for some } \Psi_{m} \in 
\sect(\hh;\O_{\hh}\twist{m})\}.
$$
This implies Theorem~\ref{th:Rreal}~(ii), by the following 
considerations.  Take a system of homogeneous coordinates 
$[\xi]=[\xi_{0},\xi']$ in $\Rb$ such that $\xio = [1,0,\dots,0]$, 
$\Rh$ is given by the equation $\xi_{0}=0$, $[\xi']$ are homogeneous 
coordinates in $\Rh$, and $q([\xi])=[\xi']$.  Then
$$
\tilde c (\psi)_{m} (\xi') = \int_{-\infty}^{+\infty} 
\psi(\sigma\xi_{0}+\xi') \sigma^m \, d\sigma .
$$
Moreover, $\sect(\hh;\O_{\hh}\twist{m})$ is precisely the space of 
homogeneous polynomials of degree $m$ in $\xi'$, considered as global 
sections of $\shc^\infty_{\Rh}\twist{m}$.

\subsection*{Notes}

{\bf \S\ref{se:RadonD}} Theorem~\ref{th:kappa-p-1} was obtained
in~\cite{D'Agnolo-Schapira96b} for $p=n-1$, where one may also find
another proof based on the Cauchy-Fantappi{\`e} formula.  The case
$p<n-1$ is from~\cite{D'Agnolo-Marastoni98a}.  In these two papers one
also finds a discussion of the case $\shm=\DP\twist{m}$, where $\shn =
\shm \dcirc\BF$ is associated to the higher dimensional analog of the
zero-rest-mass field equations.  Note that Tanisaki~\cite{Tanisaki98}
has an alternative representation theoretical proof for the case
$m=-p-1$, which extends to other flag manifolds.  This should also be
related to the work of Oshima~\cite{Oshima96}.

\medskip\noindent
{\bf \S\S\ref{ss:complexproj2}--\ref{ss:realproj}} The proofs in these 
sections are from~\cite{D'Agnolo-Schapira96b} for the case $p=n-1$, 
and from~\cite{D'Agnolo-Schapira96a} for the case $p=1$, $n=3$.  The 
general case $p<n-1$ was later obtained 
in~\cite{D'Agnolo-Marastoni98b}.

\noindent In a classical framework, Theorem~\ref{th:Rproj} can be
found in~\cite{Gelfand-Gindikin-Graev82} or \cite{Helgason84}.

\noindent The idea of investigating the real Radon transform through
the complex one, which is intrinsic in our approach, also appears
in~\cite{Gindikin98}.  For the case $p=1$, $n=3$, another close
approach can be found in~\cite{Eastwood97}.

\medskip\noindent {\bf \S\ref{ss:realaffine2}} The geometric approach
to the Cavalieri condition is from~\cite{D'Agnolo98}.  

\noindent The case $p<n-1$ is treated in~\cite{D'Agnolo-Marastoni98b}.


\begin{thebibliography}{BEW82}

\bibitem[BEW82]{Bailey-Ehrenpreis-Wells82}
  T.~Bailey, L.~Ehrenpreis, and R.~O. Wells, Jr., \emph{Weak solutions 
  of the massless field equations}, Proc.  Roy.  Soc.  London Ser.  A 
  \textbf{384} (1982), no.~1787, 403--425, \MR{84d:81021}.

\bibitem[BE89]{Baston-Eastwood89}
  R.~J. Baston and M.~G. Eastwood, \emph{The {P}enrose transform}, 
  Oxford Mathematical Monographs, The Clarendon Press Oxford 
  University Press, New York, 1989, \MR{92j:32112}.

\bibitem[Bjo93]{Bjork93}
  J.~E. Bj{\"o}rk, \emph{Analytic {$\mathcal{D}$}-modules and 
  applications}, Mathematics and its Applications, vol.~247, Kluwer 
  Academic Publishers Group, Dordrecht, 1993, \MR{95f:32014}.
  
\bibitem[Bry86]{Brylinski86}
  J.-L. Brylinski, \emph{Transformations canoniques, dualit\'e 
  projective, th\'eorie de {L}efschetz, transformations de {F}ourier 
  et sommes trigonom\'etriques}, Ast\'erisque (1986), no.~140-141, 
  3--134, 251, \MR{88j:32013}.

\bibitem[Dag98]{D'Agnolo98}
  A.~D'Agnolo, \emph{{R}adon transform and the {C}avalieri condition: 
  a cohomological approach}, Duke Math.  J. \textbf{93} (1998), no.~3, 
  597--632, \MR{1626656}.

\bibitem[DM99a]{D'Agnolo-Marastoni98a}
  A.~D'Agnolo and C.~Marastoni, \emph{Quantization of the
  {R}adon-{P}enrose transform}, preprint, 1999.

\bibitem[DM99b]{D'Agnolo-Marastoni98b}
  A.~D'Agnolo and C.~Marastoni, \emph{Real forms of the
  {R}adon-{P}enrose transform}, preprint  \textbf{218}
 Institut de Mat\'ematiques de
  Jussieu, Paris (1999).

\bibitem[DS96a]{D'Agnolo-Schapira96a}
  A.~D'Agnolo and P.~Schapira, \emph{Radon-{P}enrose transform for 
  {$\mathcal{D}$}-modules}, J. Funct.  Anal.  \textbf{139} (1996), 
  no.~2, 349--382, \MR{97h:32048}.

\bibitem[DS96b]{D'Agnolo-Schapira96b}
  A.~D'Agnolo and P.~Schapira, \emph{Leray's quantization of 
  projective duality}, Duke Math.  J. \textbf{84} (1996), no.~2, 
  453--496, \MR{98c:32012}.

\bibitem[DS98]{D'Agnolo-Schapira98}
  A.~D'Agnolo and P.~Schapira, \emph{The {R}adon-{P}enrose 
  correspondence.  {I}{I}.  {L}ine bundles and simple 
  {$\mathcal{D}$}-modules}, J. Funct.  Anal.  \textbf{153} (1998), 
  no.~2, 343--356, \MR{1614594}.

\bibitem[Eas97]{Eastwood97}
  M.~G. Eastwood, \emph{Complex methods in real integral geometry} 
  (with the collaboration of T. N.\ Bailey and C. R.\ Graham), 
  Proceedings of the 16th Winter School ``Geometry and Physics'' 
  (Srn\'\i, 1996), Rend.  Circ.  Mat.  Palermo (2) Suppl.  \textbf{46} 
  (1997), no.~3, 305--351, \MR{1469021}.

\bibitem[EPW81]{Eastwood-Penrose-Wells81}
  M.~G. Eastwood, R.~Penrose, and R.~O. Wells, Jr., \emph{Cohomology 
  and massless fields}, Comm.  Math.  Phys.  \textbf{78} (1980/81), 
  55--71, \MR{83d:81052}.

\bibitem[GGG82]{Gelfand-Gindikin-Graev82}
  I.~M. Gelfand, S.~G. Gindikin, and M.~I. Graev, \emph{Integral 
  geometry in affine and projective spaces}, Journal of Soviet Math.  
  \textbf{18} (1982), 39--167; for the Russian original see 
  \MR{82m:43017}.

\bibitem[Gin98]{Gindikin98}
  S.~G. Gindikin, \emph{Real integral geometry and complex analysis},
  Integral geometry, Radon transforms and complex analysis (Venice,
  1996), Springer, Berlin, 1998, pp.~70--98.  Lecture Notes in Math.,
  Vol.  1684, \MR{1635612}.

\bibitem[GH78]{Gindikin-Henkin78}
  S.~G. Gindikin and G.~M. Henkin, \emph{Integral geometry for 
  {$\bar\partial$}-cohomology in $q$-linearly concave domains in 
  {$\mathbb{CP}\sp{n}$}}, Funktsional.  Anal.  i Prilozhen.  
  \textbf{12} (1978), no.~4, 6--23, \MR{80a:32014}.

\bibitem[Gnc97]{Goncharov97}
  A.~B. Goncharov, \emph{Differential equations and integral 
  geometry}, Adv.  Math.  \textbf{131} (1997), no.~2, 279--343, 
  \MR{1483971}.

\bibitem[Gnz91]{Gonzalez91}
  F.~B. Gonzalez, \emph{On the range of the {R}adon $d$-plane 
  transform and its dual}, Trans.  Amer.  Math.  Soc.  \textbf{327} 
  (1991), no.~2, 601--619, \MR{92a:44002}.

\bibitem[Gri85]{Grinberg85}
  E.~L. Grinberg, \emph{On images of Radon transforms}, Duke Math.  J. 
  \textbf{52} (1985), no.~4, 939--972, \MR{87e:22020}.
  
\bibitem[GS79]{Guillemin-Sternberg79}
  V.~Guillemin and S.~Sternberg, \emph{Some problems in integral 
  geometry and some related problems in microlocal analysis}, Amer.  
  J. Math.  \textbf{101} (1979), no.~4, 915--955, \MR{82b:58087}.

\bibitem[Hel84]{Helgason84}
  S.~Helgason, \emph{Groups and geometric analysis.  Integral 
  geometry, invariant differential operators, and spherical 
  functions}, Pure and Applied Mathematics, vol.~113, Academic Press, 
  Inc., Orlando, Fla., 1984, \MR{86c:22017}.

\bibitem[HP87]{Henkin-Polyakov87}
  G.~M. Henkin and P.~L. Polyakov, \emph{Homotopy formulas for the 
  {$\overline\partial$}-operator on {$\mathbb{CP}^n$} and the 
  {R}adon-{P}enrose transform}, Math.  USSR Izvestiya \textbf{3} 
  (1987), 555--587.

\bibitem[Joh38]{John38}
  F.~John, \emph{The ultrahyperbolic differential equation with four 
  independent variables}, Duke Math.  J. \textbf{4} (1938), 300--322; 
  reprinted in {\sl 75 years of {R}adon transform} ({V}ienna, 1992), 
  {C}onf.  {P}roc.  {L}ecture {N}otes {M}ath.  {P}hys., {IV}, 
  {I}nternat.  {P}ress, {C}ambridge, {M}A, 1994, \MR{96d:35096}.

\bibitem[Kak97]{Kakehi97}
  T.~Kakehi, \emph{Range theorems and inversion formulas for {R}adon 
  transforms on {G}rassmann manifolds}, Proc.  Japan Acad.  Ser.  A 
  Math.  Sci.  \textbf{73} (1997), no.~5, 89--92, \MR{1470177}.

\bibitem[Kas84]{Kashiwara84}
  M.~Kashiwara, \emph{The {R}iemann-{H}ilbert problem for holonomic 
  systems}, Publ.  Res.  Inst.  Math.  Sci.  \textbf{20} (1984), 
  no.~2, 319--365, \MR{86j:58142}.

\bibitem[Kas70]{Kashiwara95}
  M.~Kashiwara, \emph{Algebraic study of systems of partial 
  differential equations} (a translation by A.~D'Agnolo and 
  J.-P.~Schneiders of Kashiwara's Master's Thesis, Tokyo University, 
  December 1970), M\'em.  Soc.  Math.  France (N.S.) (1995), no.~63, 
  xiv+72, \MR{97f:32012}.

\bibitem[KSa90]{Kashiwara-Schapira90}
  M.~Kashiwara and P.~Schapira, \emph{Sheaves on manifolds}, 
  Grundlehren der Mathematischen Wissenschaften, vol.  292, 
  Springer-Verlag, Berlin, 1990, \MR{92a:58132}.

\bibitem[KSa96]{Kashiwara-Schapira96}
  M.~Kashiwara and P.~Schapira, \emph{Moderate and formal cohomology 
  associated with constructible sheaves}, M\'em.  Soc.  Math.  France 
  (N.S.) (1996), no.~64, iv+76, \MR{97m:32052}.

\bibitem[KSm94]{Kashiwara-Schmid94}
  M.~Kashiwara and W.~Schmid, \emph{Quasi-equivariant 
  {$\mathcal{D}$}-modules, equivariant derived category, and 
  representations of reductive {L}ie groups}, Lie theory and geometry, 
  Progr.  Math., vol.  123, Birkh\"auser Boston, Boston, MA, 1994, 
  pp.~457--488, \MR{96e:22031}.

\bibitem[Mar67]{Martineau67}
  A.~Martineau, \emph{Equations diff{\'e}rentielles d'ordre infini},
  Bull.  Soc.  Math.  France \textbf{95} (1967), 109--154.

\bibitem[Osh96]{Oshima96}
  T.~Oshima, \emph{Generalized {C}apelli identities and boundary value 
  problems for {$\mathrm{GL}(n)$}}, Structure of solutions of 
  differential equations (Katata/Kyoto, 1995), World Sci.  Publishing, 
  River Edge, NJ, 1996, pp.~307--335, \MR{98f:22021}.

\bibitem[Rad17]{Radon17}
  J.~Radon, \emph{{\"U}ber die Bestimmung von Funktionen durch ihre 
  Integralwerte l{\"a}ngs gewisser Mannigfaltigkeiten.}, Ber.  Verh.  
  Konigl.  Sachs.  Ges.  Wiss.  Leipzig \textbf{69} (1917), 262 -- 
  277; reprinted in {\sl 75 years of {R}adon transform} ({V}ienna, 
  1992), {C}onf.  {P}roc.  {L}ecture {N}otes {M}ath.  {P}hys., {IV}, 
  {I}nternat.  {P}ress, {C}ambridge, {M}A, 1994, \MR{95m:44007}.

\bibitem[SKK73]{Sato-Kawai-Kashiwara73}
  M.~Sato, T.~Kawai, and M.~Kashiwara, \emph{Microfunctions and 
  pseudo-differential equations}, Hyperfunctions and 
  pseudo-differential equations (Proc.  Conf., Katata, 1971; dedicated 
  to the memory of Andr\'e Martineau), Springer, Berlin, 1973, 
  pp.~265--529.  Lecture Notes in Math., Vol.  287, \MR{54:8747}.

\bibitem[Scp85]{Schapira85}
  P.~Schapira, \emph{Microdifferential systems in the complex domain}, 
  Grundlehren der mathematischen Wissenschaften, vol.~269, 
  Springer-Verlag, Berlin-New York, 1985, \MR{87k:58251}.
  
\bibitem[SS94]{Schapira-Schneiders94}
  P.~Schapira and J.-P.~Schneiders, \emph{Index theorem for elliptic 
  pairs}, Ast{\'e}risque, vol.~224, 1994, \MR{96a:58179}.
  
\bibitem[Scn94]{Schneiders94}
  J.-P. Schneiders, \emph{An introduction to {$\mathcal{D}$}-modules}, 
  Algebraic Analysis Meeting (Li{\`e}ge, 1993), Bull.  Soc.  Roy.  
  Sci.  Li{\`e}ge \textbf{63} (1994), no.~3-4, 223--295, 
  \MR{95m:32019}.

\bibitem[Sek96]{Sekiguchi96}
  H.~Sekiguchi, \emph{The {P}enrose transform for certain non-compact 
  homogeneous manifolds of {${\rm {U}}(n,n)$}}, J. Math.  Sci.  Univ.  
  Tokyo \textbf{3} (1996), no.~3, 655--697, \MR{98d:22011}.
  
\bibitem[Tan98]{Tanisaki98}
  T.~Tanisaki, \emph{Hypergeometric systems and {R}adon transforms for 
  hermitian symmetric spaces}, to appear in Adv. Studies in Pure 
  Math., Vol.  26.

\bibitem[Wel81]{Wells81}
  R.~O. Wells, Jr., \emph{Hyperfunction solutions of the 
  zero-rest-mass field equations}, Comm.  Math.  Phys.  \textbf{78} 
  (1980/81), no.~4, 567--600, \MR{83d:81053}.

\end{thebibliography}
\end{document}